\documentclass[11pt]{article}
\usepackage{amsfonts,amsthm}
\usepackage{eucal}
\usepackage{amssymb,epsfig,latexsym}
\usepackage{graphicx,amsmath,amssymb,latexsym,psfrag}
\usepackage{graphics,graphicx}
\usepackage[colorlinks=true, urlcolor=blue,linkcolor=blue, citecolor=blue]{hyperref}

\newtheorem{Theorem}{Theorem}[section]
\newtheorem{Proposition}[Theorem]{Proposition}

\newtheorem{Definition}[Theorem]{Definition}
\newtheorem{Corollary}[Theorem]{Corollary}
\newtheorem{Lemma}[Theorem]{Lemma}
\newtheorem{Remark}[Theorem]{Remark}
\newtheorem{Example}[Theorem]{Example}
\def\ind{{\rm 1\hspace{-0.90ex}1}}

\newcommand{\R}{{\mathbb R}}

\newcommand{\E}{{\mathbb E}}

\renewcommand{\d}{{\operatorname d}}

\newcommand{\dM}{{\mathrm d}}

\renewcommand{\d}{{\operatorname d}}

\newcommand{\card}{\mathrm{card}}

\def\ind{{\rm 1\hspace{-0.90ex}1}}

\title{\huge Gaussian Approximation and Moderate Deviations of Poisson Shot Noises with Application to Compound Generalized Hawkes Processes}
\author{
Mahmoud Khabou\thanks{Imperial College London.  M.K. was supported by the Project EDDA (ANR-20-IADJ-0003) of the French National Research Agency (ANR)
e-mail: \tt{m.khabou@imperial.ac.uk}}
\,\,and\,\,
Giovanni
Luca Torrisi\thanks{Consiglio Nazionale delle Ricerche
(CNR),  Via dei Taurini 19, 00185 Roma,  Italy. G.L.T was supported by group GNAMPA of INdAM.
e-mail: \tt{giovanniluca.torrisi@cnr.it}}}
\allowdisplaybreaks

\begin{document}

\maketitle

\begin{abstract}

In this article, we give explicit bounds on the Wasserstein and the Kolmogorov distances between random variables lying in the first chaos of the Poisson space and the standard Normal distribution, using the results proved in \cite{LPS}. Relying on the theory developed in \cite{SS} and on a fine control of the cumulants of the first chaoses,  we also derive moderate deviation principles,  Bernstein-type concentration inequalities and Normal approximation bounds with Cram\'er correction terms for the same variables. The aforementioned results are then applied to Poisson shot-noise processes and,  in particular,  to
the generalized compound Hawkes point processes (a class of stochastic models, introduced in this paper,  which generalizes classical Hawkes processes).  This extends the recent results in \cite{HHKR,KPR} regarding the Normal approximation and in \cite{ZHU} for moderate deviations. 

%We provide explicit bounds on the Wasserstein and the Kolmogorov distances between generalized compound Hawkes point processes (a class of stochastic models, introduced in this paper,  which generalizes classical Hawkes processes) and the standard Normal law,  extending the recent results in \cite{HHKR,KPR}.  For generalized compound Hawkes point processes we provide also moderate deviations,  Bernstein-type concentration inequalities and Normal approximation bounds with Cram\'er correction term,  extending the result in \cite{ZHU}.  Our approach exploits the Normal approximation bounds for functionals of the Poisson measure provided in \cite{LPS} and,  as far as moderate deviations is concerned,  is based on a suitable application of the theory developed in \cite{SS} via a fine control of the cumulants of the first chaos on the Poisson space.
\end{abstract}

\section{Introduction}

Shot noise models built on one-dimensional Poisson processes are very popular in applied probability.  Due to their versatility and mathematical tractability,  they find applications in many fields,  such as insurance,  finance,  queueing theory and neuroscience (see e.g.  \cite{B, BD, GMT1, GMT2, KL,  KM0,KM, LT, MT, MST, P, TL1, T}).  Shot noise models, whose underlying point processes are spatial Poisson processes (see Section \ref{sec:PSNmodel} for a formal definition) are a bit less popular,  but they play an important role in wireless communication,  where they are exploited as models of the inference in ad hoc networks (see e.g.  \cite{BB, BB1, BB2, GT2,PT,TL}). Furthermore,  as explained in detail in the next section,  shot noise models,  whose underlying point processes are spatial Poisson processes (hereafter also called \lq\lq spatial\rq\rq\, Poisson shot noise models) encompass spatial Poisson cluster point processes,  which are widely used in many research areas such as spatial statistics (see e.g.  \cite{MW}).  Since \lq\lq spatial\rq\rq\,
Poisson shot noise models are stochastic integrals with respect to a Poisson random measure,  Gaussian approximation bounds for the Wasserstein and the Kolmogorov distances between such random variables (properly standardized) and the standard Normal law can be easily obtained applying the general theory developed in the seminal paper \cite{LPS}.  One of the main achievements of the present article are explicit bounds for the Wasserstein and Kolmogorov distances between a properly standardized compound sum,
which extends Poisson cluster and Hawkes point processes,  and the standard Normal law (see Corollaries \ref{cor:2},  \ref{cor:11102023quarto} and \ref{cor:11102023quartoBIS}).
These results improve and go beyond the findings in \cite{HHKR,KPR},  exploiting a considerably simpler approach (see also the discussion in Section \ref{subsec:21122023primo}).  Using a well-known link between cumulants and large deviations theory (see \cite{SS}),  we also provide sufficient conditions which guarantee moderate deviations,  Bernstein-type concentration inequalities and Normal approximation bounds with Cram\'er correction term (see Definition \ref{def:27102023} for details) for sequences of random variables which belong to the first chaos on the Poisson space (see Theorem \ref{thm:SchulteThale}).  Then,  we transfer such results to sequences of 
\lq\lq spatial\rq\rq\, Poisson shot noise models.  As a main application,  we provide moderate deviations,  Bernstein-type concentration inequalities and Normal approximation bounds with Cram\'er correction term for sequences of compound sums,  which extend Poisson cluster and Hawkes point processes (see Corollaries \ref{cor:ST},  \ref{cor:11102023quinto} and \ref{cor:11102023nono}). Remarkably,  the result on moderate deviations recovers, under an alternate condition on the fertility function,  the moderate deviations for the number of points of a classical Hawkes process on the time interval $(0,t]$ proved in \cite{ZHU} (see Section \ref{subsec:21122023secondo}).

The paper is structured as follows. In Section \ref{sec:notations} we introduce the Poisson shot noise models considered in the paper,  and show that compound Poisson cluster point processes
and generalized compound Hawkes processes are indeed particular Poisson shot noise models.  Furthermore,  we recall a simple model of wireless communication, which accounts for interference effects described by a Poisson shot noise.  In Section \ref{sec:Informal} we provide an informal description of our results.  In Section \ref{sec:prel} we give Gaussian approximation bounds for the Wasserstein and the Kolmogorov distances between a random variable belonging to the first chaos of the Poisson space and the standard Normal law,  and we give moderate deviations,  Bernstein-type concentration inequalities and Normal approximation bounds with Cram\'er correction term for sequences of random variables belonging to the first chaos on the Poisson space.  Applications of the results of Section \ref{sec:prel} to \lq\lq spatial\rq\rq\, Poisson shot noise models and compound Poisson cluster point processes are provided in Sections \ref{sec:APPPSN} and \ref{sec:CPCP},  respectively.   The general results on Gaussian approximation and moderate deviations are applied to generalized compound Hawkes processes in Sections \ref{sec:CHP} and \ref{sec:WIRE},  respectively.

\section{Poisson shot noise random variables}\label{sec:PSNmodel}
\label{sec:notations}
 Throughout this article, if $x$ is a point in some set $E$ and $C \subset E$, then $C-x$ denotes the set
$\{y-x, y \in C\}.$
A Poisson shot noise random variable is a real-valued random variable of the form
\begin{equation}\label{eq:PSN}
S(C):=\sum_{n\geq 1}H(C-X_n,Z_n),\quad\text{$C\in\mathcal{B}(\mathbb R^d)$.}
\end{equation}
Here $\mathcal{B}(\R^d)$ denotes the Borel $\sigma$-field on $\R^d$,  $d\geq 1$,  $\mathcal{P}\equiv\{(X_n,Z_n)\}_{n\geq 1}$
is a Poisson process on $\mathbb R^d\times\bold Z$ with mean measure $\lambda(x)\mathrm{d}x\mathbb{Q}(\mathrm{d}z)$, 
$(\bold Z,\mathcal Z)$ is a measurable space,
$\lambda:\mathbb R^d\to [0,\infty)$ is a locally integrable intensity function,  $\mathbb{Q}$ is a probability measure on $\bold Z$
and $H:\mathcal{B}(\R^d)\times\bold Z\to\R$ is a mapping such that,
for each fixed $C\in\mathcal{B}(\mathbb R^d)$,  the function
\[
(x,z)\in\mathbb R^d\times\bold{Z}\mapsto H(C-x,z)\in\mathbb R
\]
is measurable.  Poisson shot noise random variables encompass a variety of important stochastic models.

\subsection{Compound Poisson cluster point processes}\label{sec:markedPC}

Let $\{X_n\}_{n\geq 1}$ be the points of a Poisson process on $\mathbb R^d$,  $d\geq 1$,  with a locally integrable intensity function $\lambda:\mathbb{R}^d\to [0,\infty)$ and let $\{Z_n(\cdot,\cdot)\}_{n\geq 1}$ be a sequence of independent and identically distributed simple point processes on $\R^d \times \R$,  
independent of $\{X_n\}_{n\geq 1}$. More concretely, for $(C_1,C_2) \in \mathcal B (\R^d) \times \mathcal B (\R)$, $Z_n(C_1,C_2)$ counts the number of points of the $n$-th point process that fall in $C_1$ and whose marks are in $C_2$. For each $n\geq 1$,  we denote the points of
$Z_n(\cdot, \cdot)$ by $\{(Y_{n,k},M_{n,k})\}_{k\geq 0}$,  and we assume that
$Y_{n,0}:= \mathbf {0}$ (which implies $Z_n(\{\bold 0\},\R):=1$) and that the sequence $\{M_{n,k}\}_{k\geq 0}$
is independent of $\{Y_{n,k}\}_{k\geq 0}$.  Furthermore,  we suppose that
the random variables $\{M_{n,k}\}_{n\geq 1,\,k\geq 0}$ are independent and identically distributed.
Throughout this paper we denote by $M$ the generic random variable $M_{n,k}$.
  
One naturally interprets the first components of the points of $Z_n(\cdot,\cdot)$ as \lq\lq locations\rq\rq\, and the second components as \lq\lq marks\rq\rq,\, which describe some characteristic of the location to which is \lq\lq attached\rq\rq.  Hereafter,  for $n\geq 1$,  we consider the point processes 
$\theta_{X_n}Z_n(\cdot,\cdot)\equiv\{(X_n+Y_{n,k},M_{n,k})\}_{k\geq 0}$.

For arbitrarily fixed $n\geq 1$ and $C \in \mathcal B (\R^d)$,  we define the random variable 
\begin{equation}\label{def:upsilon}
\upsilon(Z_n)(C):=\sum _{k=0}^{Z_n(C,\R)-1}M_{n,k},
\end{equation}
which aggregates the marks attached to the locations that fall in $C$.  It turns out that the random variable, say $V(C)$, which aggregates all
the marks attached to the points,  which fall in $C$,  of the Poisson cluster point process 
\[
N\equiv\bigcup_{n\geq 1}\{X_n+Y_{n,k}\}_{k\geq 0}
\]
is a Poisson shot noise random variable.  Indeed,
\begin{equation}\label{eq:rapprS}
V(C):=\sum_{n\geq 1}\sum _{k=0}^{\theta_{X_n}Z_n(C,\R)-1}M_{n,k}=
\sum_{n\geq 1}\upsilon(\theta_{X_n}Z_n)(C)=\sum_{n\geq 1}\upsilon(Z_n)(C-X_n),
\end{equation}
is a random variable of the form \eqref{eq:PSN} with $H(C-x,z):=v(z)(C-x)$,  $x\in\mathbb{R}^d$,  $z\in\bold{Z}:=\bold{N}_{\mathbb{R}^d\times\mathbb{R}}$.
Here $\bold{N}_{\mathbb R^d\times\mathbb R}$ denotes the space of $\sigma$-finite counting measures on $(\mathbb R^d\times\mathbb{R},\mathcal B(\mathbb R^d)\otimes\mathcal B(\mathbb R))$ equipped with the 
usual $\sigma$-field (see Section \ref{sec:prel} for details),  and
\begin{equation}\label{eq:15122023pom5}
v(z)(C):=\sum_{k\geq 0}\bold{1}_C(y_k)m_k,\quad\text{for $z\equiv\{(y_k,m_k)\}_{k\geq 0}$.} 
\end{equation}
Note that if $M_{n,k}:=1$ for every $n\geq 1$ and $k\geq 0$,  then the random variable
\begin{equation}\label{eq:13122023}
N(C):=V(C)= \sum_{n\geq 1}Z_n(C-X_n,\R)
\end{equation}
equals the number of points of the Poisson cluster point process $N$ which fall in $C\in\mathcal{B}(\R^d)$.

\subsection{Generalized Hawkes processes and generalized compound Hawkes processes}

Let $N\equiv\{N(C)\}_{C\in\mathcal{B}(\R^d)}$ be the Poisson cluster point process defined by \eqref{eq:13122023}.  $N$ will be named {\em generalized Hawkes process}
if the random variable $Z:=Z_1(\R^d,\R)$ is distributed as the total progeny of a sub-critical Galton-Watson process with one ancestor. 

It is worthwhile to note that:\\
\noindent $(i)$ {\em Classical Hawkes processes} on $(0,\infty)$ (respectively,  on  $\mathbb R$) with parameters $(\lambda,g)$,  introduced in the seminal papers \cite{H,HO},  are particular generalized Hawkes processes.  Indeed,  they are Poisson cluster point processes defined as follows: $(1)$ The process of cluster centers $\{X_n\}_{n\geq 1}$ is a Poisson process on $(0,\infty)$ (respectively,  on $\mathbb R$) with constant intensity equal to $\lambda>0$; $(2)$ The points of the cluster $\theta_{X_n}Z_n(\cdot,\R)$
are partitioned into generations and generated recursively as follows. 
The ancestor constitutes the generation $0$ of the cluster and it is located at $X_n$.
Given $X_n$,  the ancestor generates points of the $1$st generation of the cluster according to a non-homogeneous Poisson process on $(X_n,\infty)$ with intensity function $g(\cdot-X_n)$, where $g:\mathbb R\to [0,\infty)$ is a measurable function which is null on $(-\infty,0]$ and such that $h:=\int_0^\infty g(x)\mathrm{d}x<1$.  In turn,  given the points of the $1$st generation of the cluster,  a point of this generation, which is located at $X$,  generates points of the second generation of the cluster according to a non-homogeneous Poisson process on $(X,\infty)$ with intensity function $g(\cdot-X)$.  And so on and so forth.  Note that $Z_n(\mathbb R,\R)=\theta_{X_n}Z_n(\mathbb R,\R)=\theta_{X_n}Z_n([X_n,\infty),\R)=Z_n([0,\infty),\R)$ is distributed as the total progeny of a sub-critical Galton-Watson process with one ancestor and Poisson offspring law with mean $h$;\\
\noindent $(ii)$ {\em Spatial Hawkes processes on $\R^d$, } $d\geq 1$, with parameters $(\lambda,g)$,  introduced in \cite{BMR} and further studied in \cite{MT},  are particular generalized Hawkes processes.  Indeed,  they are Poisson cluster point processes defined as follows:  $(1)$ The process of cluster centers $\{X_n\}_{n\geq 1}$ is a Poisson process on $\mathbb R^d$,  $d\geq 1$,  with constant intensity equal to $\lambda>0$; $(2)$ The points of the cluster $\theta_{X_n}Z_n(\cdot,\R)$ are partitioned into generations and generated recursively as follows. 
The ancestor constitutes the generation $0$ of the cluster and it is located at $X_n$.
Given $X_n$,  the ancestor generates points of the $1$st generation of the cluster according to a non-homogeneous Poisson process on $\mathbb R^d$ with intensity function $g(\cdot-X_n)$, where $g:\mathbb R^d\to [0,\infty)$ is a measurable function such that $h:=\int_{\mathbb R^d}g(x)\mathrm{d}x<1$.  In turn,  given the points of the $1$st generation of the cluster,  a point of this generation,  which is located at $X$,  generates points of the second generation of the cluster according to a non-homogeneous Poisson process on $\mathbb R^d$ with intensity function $g(\cdot-X)$.  And so on and so forth.  Note that $\theta_{X_n}Z_n(\mathbb R^d,\R)=Z_n(\mathbb R^d,\R)$ is distributed as the total progeny of a sub-critical Galton-Watson process with one ancestor and Poisson offspring law with mean $h$. 

Note that,  according to these definitions,  classical Hawkes processes on $\mathbb R$
are different from \lq\lq spatial\rq\rq\, Hawkes processes on $\mathbb R$.

The collection of random variables $V\equiv\{V(C)\}_{C\in\mathcal{B}(\R^d)}$, where $V(C)$ is defined by \eqref{eq:rapprS},  will be named {\em generalized compound Hawkes process} if the random variable $Z$ is distributed as the total progeny of a sub-critical Galton-Watson process with one ancestor.  Note that $V(C)$ aggregates the marks attached to the points of a generalized Hawkes process which fall in $C$.

\subsection{Interference in wireless communication}\label{sec:07122023}

Consider the following simple model of wireless communication,  which accounts for interference effects that arise when several nodes transmit at the same time.
Suppose that transmitting nodes (e.g.,  antennas) are located according to $\{X_n\}_{n\geq 1}$,  a Poisson process on the plane with intensity function $\lambda(\cdot)$,  i.e.,  $X_n$ is the location of node $n$,  and denote by $Z_n\in (0,\infty)$ the signal power of the transmitting node $n$.  Suppose that the sequence $\{Z_n\}_{n\geq 1}$ is independent of the Poisson process,  and that the random variables $Z_n$,  $n\geq 1$,  are independent and identically distributed.  Assume that a receiver is located at the origin
$\bold 0\in\mathbb R^2$ and that a new transmitter is added at $x\in\mathbb R^2$ and has signal power $y\in (0,\infty)$.

Suppose that the physical propagation of the signal is described by a measurable positive function $A:\mathbb R^2\to (0,\infty)$,  which gives the attenuation or path loss of the signal power.  For simplicity,  we assume that the random fading (due to occluding objects,  reflections,  multipath interference,  etc. ...) is encoded in the random variables $Z_n\in (0,\infty)$.  Thus,  $Z_n A(X_n)$ is the received power at the origin due to the transmitting node at  $X_n$,  and the total interference at the origin,  due to simultaneous transmissions,  is equal to 
\[
I(\{\bold 0\}):=\sum_{n\geq 1}Z_n A(X_n).
\]
Note that this is a Poisson shot noise random variable of the form \eqref{eq:PSN}.  Indeed,  let $H:\mathcal{B}(\R^2)\times (0,\infty)\to (0,\infty)$ be a mapping which,  restricted on $\R^2\times (0,\infty)$,
coincides with $\widetilde{H}(x,z):=zA(-x)$.  Then 
\[
S(\{\bold 0\})=\sum_{n\geq 1}H(\{\bold 0\}-X_n,Z_n)=\sum_{n\geq 1}\widetilde{H}(-X_n,Z_n)=I(\{\bold 0\}).
\]

We refer the reader to \cite{BB1, BB2} for more insight into this model,  and limit ourselves to say that the receiver at the origin can decode the signal of power $y\in (0,\infty)$ from the transmitter at $x\in\mathbb R^2$ if and only if the Signal to  Interference plus Noise Ratio (SINR) is bigger than a given threshold, i.e.,
\[
\mathrm{SINR}:=\frac{yA(x)}{I(\{\bold 0\})+w}\geq\tau,
\]
where $w$ is e.g. a thermal noise near the receiver at the origin and $\tau$ is the given threshold.

\section{Informal description of the results}\label{sec:Informal}

We start noticing that some results in this paper refer to sequences of Poisson shot noise random variables of the form
\begin{equation}\label{eq:PSNsequence}
S_\ell(C_\ell):=\sum_{n\geq 1}H(C_\ell-X_n^{(\ell)},Z_n^{(\ell)}),\quad\text{$\ell\geq 1$, $\{C_\ell\}_{\ell\geq 1}\subset\mathcal{B}(\mathbb R^d)$}
\end{equation}
where,  for each $\ell\geq 1$,  $\mathcal{P}_\ell=\{(X_n^{(\ell)},Z_n^{(\ell)})\}_{n\geq 1}$ is a Poisson process on $\mathbb{R}^d\times\bold Z$ with mean measure
$\lambda_\ell(x)\mathrm{d}x\mathbb{Q}_\ell(\mathrm{d}z)$,  $\lambda_\ell:\mathbb{R}^d\to [0,\infty)$ is a locally integrable function and $\mathbb{Q}_\ell(\cdot)$ is a probability measure on a measurable space $(\bold Z,\mathcal{Z})$.

The achievements of the paper are the following.\\ 
\noindent $(i)$ We provide bounds for the Wasserstein and the Kolmogorov distances (hereafter denoted by $d_W$ and $d_K$,  respectively; see Section \ref{sec:27102023} for a formal definition of these probability metrics) between a standard Gaussian random variable $G$ and
\begin{equation}\label{eq:11102023terzo}
T(C):=\frac{S(C)-\mathbb E S(C)}{\sqrt{\mathbb{V}\mathrm{ar}(S(C))}}, \quad C\in\mathcal{B}(\mathbb R^d)
\end{equation}
see Theorem \ref{thm:GaussPoissoncluster}.  As particular cases (see Remark \ref{re:particularcases}),  we get bounds for the Wasserstein and the Kolmogorov distances between $G$ 
and

\begin{equation}\label{eq:11102023terzoV}
W(C):=\frac{V(C)-\mathbb E V(C)}{\sqrt{\mathbb{V}\mathrm{ar}(V(C))}}, \quad C\in\mathcal{B}(\mathbb R^d)
\end{equation}
where $V(C)$ is defined by \eqref{eq:rapprS},  and between $G$ and
\begin{equation}\label{eq:11102023terzoI}
L(\{\bold 0\}):=\frac{I(\{\bold 0\})-\mathbb E I(\{\bold 0\})}{\sqrt{\mathbb{V}\mathrm{ar}(I(\{\bold 0\}))}}, \quad\bold 0\in\mathbb R^2
\end{equation}
where $I(\{\bold 0\})$ is defined in Section \ref{sec:07122023}.\\
\noindent$(ii)$ We provide moderate deviations,  Bernstein-type concentration inequalities and Normal approximation bounds with Cram\'er correction term for the sequence 
$\{T_\ell(C_\ell)\}_{\ell\geq 1}$, where
\[
T_\ell(C_\ell):=\frac{S_\ell(C_\ell)-\mathbb E S_\ell(C_\ell)}{\sqrt{\mathbb{V}\mathrm{ar}(S_\ell(C_\ell))}}, \quad C_\ell\in\mathcal{B}(\mathbb R^d)
\]
see Theorem \ref{thm:ModeratePoissoncluster}.  As particular cases (see Remark \ref{re:particularcasesbis}),  we get moderate deviations,  Bernstein-type concentration inequalities and Normal approximation bounds with Cram\'er correction term for the sequences:

\begin{equation}\label{eq:11102023terzoVelle}
W_\ell(C_\ell):=\frac{V_\ell(C_\ell)-\mathbb E V_\ell(C_\ell)}{\sqrt{\mathbb{V}\mathrm{ar}(V_\ell(C_\ell))}}, \quad C_\ell\in\mathcal{B}(\mathbb R^d),\quad\ell\geq 1
\end{equation}
where $V_\ell(C_\ell)$ is defined by \eqref{eq:rapprS},  with obvious modifications,  and
\begin{equation}\label{eq:11102023terzoIelle}
L_\ell(\{\bold 0\}):=\frac{I_\ell(\{\bold 0\})-\mathbb E I_\ell(\{\bold 0\})}{\sqrt{\mathbb{V}\mathrm{ar}(I_\ell(\{\bold 0\}))}}, \quad\bold 0\in\mathbb R^2
\end{equation}
where $I_\ell(\{\bold 0\})$ is defined in Section \ref{sec:07122023}, with obvious modifications.\\
\noindent $(iii)$ If the Poisson process $\{X_n\}_{n\geq 1}$ has intensity function of the form
$\lambda(x):=\lambda\bold{1}_{B}(x)$,  $x\in\mathbb R^d$,  for some positive constant $\lambda>0$ and a suitable Borel set $B\subseteq\mathbb R^d$,  then the bounds on $d_W(W(C),G)$ and $d_K(W(C),G)$ are particularly simple and depend only on $\lambda$,  the Lebesgue measure of $B\cap C$ and a few moments of $M$ and  $Z$,  see Corollary \ref{cor:2}.  If $Z$ is distributed as the total progeny of a sub-critical Galton-Watson process with one ancestor,
then we are able to compute the moments of $Z$ in terms of the moments of the offspring law,  see Proposition \ref{thm:branching}. This allows for explicit bounds when
$V=\{V(C)\}_{C\in\mathcal{B}(\R^d)}$ is a generalized compound Hawkes process with Poisson or Binomial offspring laws,  see Corollaries
\ref{cor:11102023quarto} and \ref{cor:11102023quartoBIS},  respectively.\\
\noindent $(iv)$ If the Poisson process $\{X_n^{(\ell)}\}_{n\geq 1}$ has intensity function of the form
$\lambda_\ell(x)=\lambda_\ell \boldsymbol 1_{B_\ell}(x)$,  $x\in\mathbb R^d$, for  
positive constants $\lambda_\ell>0$ and suitable Borel sets
$B_\ell\in\mathcal{B}(\mathbb R^d)$, and 
$\mathbb Q_\ell\equiv\mathbb Q$,  for each $\ell\geq 1$,  then the condition which guarantees a moderate deviation principle,  a Bernstein-type concentration inequality and a Normal approximation bound with Cram\'er correction term for the sequence $\{W_\ell(C_\ell)\}_{\ell\geq 1}$ is quite simple
(see condition \eqref{hyp:29092023} of Corollary \ref{cor:ST}),  and it allows for applications to generalized compound Hawkes processes with Poisson and Binomial offspring distribution,  see Corollaries \ref{cor:11102023quinto} and \ref{cor:11102023nono},  respectively.

We conclude this section emphasizing that the basic idea of this paper is very simple (expecially if compared with the techniques exploited in \cite{HHKR, KPR} for the Gaussian approximation of classical Hawkes processes on $(0,\infty)$).  Since the random variable $S(C)$ can be rewritten as the Poisson integral
\begin{equation}
\label{eq:Pois_integral}
S(C)=\int_{\mathbb{R}^d\times\bold Z}H(C-x,z)\mathcal{P}(\mathrm{d}x,\mathrm{d}z),\quad C\in\mathcal{B}(\R^d):
\end{equation}

\noindent$(i)$ We apply the quantitative central limit theorems,  in the Wasserstein and the Kolmogorov metrics,  for functionals of the Poisson measure proved in the seminal paper \cite{LPS},  to provide Normal approximation bounds for first chaoses on the Poisson space (see Theorems \ref{thm:LPS1}),  and then we apply such bounds to the random variable $T(C)$.\\
\noindent$(ii)$ We prove moderate deviations,  Bernstein-type concentration inequalities and Normal approximation bounds with Cram\'er correction for a sequence of first chaoses on the Poisson space (see e.g. Theorem \ref{thm:SchulteThale}), and then we apply such results to the sequence $\{T_\ell(C_\ell)\}_{\ell\geq 1}$.

\section{Gaussian approximation and moderate deviations of the first chaos on the Poisson space}\label{sec:prel}

Let $(A,\mathcal A,\alpha)$ be a measure space with $\alpha(\cdot)$ a $\sigma$-finite measure and let $\bold N_A$ be the set of all $\sigma$-finite counting measures on $(A,\mathcal A)$ equipped with the $\sigma$-field generated by the mappings
$\nu\mapsto\nu(B)$, $B\in\mathcal A$.  $\Pi$ is a Poisson measure on $(A,\mathcal A)$ with mean measure $\alpha(\cdot)$ if it is a measurable mapping from an underlying probability space
$(\Omega,\mathcal F,\mathbb P)$ to $\bold N_A$ such that: $(i)$ For any $B\in\mathcal A$,  $\Pi(A)$ is Poisson distributed with mean $\alpha(A)$; $(ii)$ If $B_1,\ldots,B_n\in\mathcal A$,
$n\in\mathbb N$,  are pairwise disjoint,  then the random variables $\Pi(B_1),\ldots,\Pi(B_n)$ are independent.

Let $\{\Pi_\ell\}_{\ell\in\mathbb N}$ be a sequence of Poisson measures on $(A,\mathcal A)$,  defined on the probability space $(\Omega,\mathcal F,\mathbb P)$.   Suppose that $\Pi_\ell$ has a $\sigma$-finite mean measure $\alpha_\ell(\cdot)$,  $\ell\in\mathbb N$.  We denote by $L^m(\alpha_\ell)$ the space of measurable functions $f:A\to\R$ such that $\int_A|f(a)|^m\alpha_\ell(\d a)<\infty$,  $m\in\mathbb N$,  and, for $\{f_\ell\}_{\ell\in\mathbb N}\in L^2(\alpha_\ell)$,  we consider
\[
I^{(\ell)}(f_\ell):=\int_A f_\ell(a)(\Pi_\ell(\d a)-\alpha_\ell(\d a)),\quad\ell\in\mathbb N
\]
the first chaos on the Poisson space,  i.e.,  the first order stochastic integral of $f_\ell$ with respect to the compensated Poisson measure $\Pi_\ell(\d a)-\alpha_\ell(\d a)$.  If the law of $\Pi_\ell$ does not depend on $\ell$,  we suppress the dependence on $\ell$ of the related quantities and,  for $f\in L^2(\alpha)$,  we simply write
\[
I(f):=\int_A f(a)(\Pi(\d a)-\alpha(\d a)).
\]

\subsection{Bounds on the Wasserstein and Kolmogorov distances between the law of a first chaos on the Poisson space and the standard Normal distribution}\label{sec:27102023}

Let $X$ and $Y$ be two real-valued random variables defined on $(\Omega,\mathcal{F},\mathbb P)$.  The Wasserstein and the Kolmogorov distances between
the law of $X$ and the law of $Y$,  written $d_W(X,Y)$ and $d_K(X,Y)$,  respectively, 
are defined by
\[
d_W(X,Y):=\sup_{g\in\mathrm{Lip}(1)}|\mathbb E[g(X)-g(Y)]|
\]
and
\[
d_K(X,Y):=\sup_{x\in\R}|\mathbb{P}(X\leq x)-\mathbb{P}(Y\leq x)|.
\]
Here, $\mathrm{Lip}(1)$ denotes the set of Lipschitz functions $g:\R\to\R$ with Lipschitz constant less than or equal to $1$.
We recall that throughout this paper,  $G$ denotes a random variable distributed according to the standard Normal law.

\begin{Theorem}\label{thm:LPS1}
Let $f\in L^2(\alpha)$ be such that $\|f\|_{L^2(\alpha)}=1$.  Then
\[
d_W(I(f),G)\leq\int_{A}|f(a)|^3\alpha(\mathrm{d}a)
\]
and
\begin{align}
&d_K(I(f),G)
\leq\left(1+\frac{1}{2}\max\left\{4,
\left[4\int_A|f(a)|^4\alpha(\mathrm{d}a)+2\right]^{1/4}
\right\}\right)\int_A|f(a)|^3\alpha(\mathrm{d}a)
+\left(\int_A|f(a)|^4\alpha(\mathrm{d}a)\right)^{1/2}.\nonumber
\end{align}
\end{Theorem}

\noindent{\it Proof.}\\
We refer the reader to \cite{LP} for all the notions of stochastic analysis on the Poisson space used in this proof.  Let $F$ be a functional of $\Pi$, i.e.,  $F=\mathfrak{f}(\Pi)$, where $\mathfrak f$ is a real-valued measurable function defined on $\bold N_A$.  
We recall that the difference operator $D$ is defined by
\[
D_a F=\mathfrak{f}(\Pi+\delta_a)-\mathfrak{f}(\Pi),\quad a\in A
\]
where $\delta_a$ is the Dirac measure at $a\in A$,  and that the second difference operator $D^2$ is defined by
\[
D_{a_1,a_2}^{2}F=D_{a_1}(D_{a_2}F),\quad a_1,a_2\in A. 
\]
We also recall that the domain of $D$,  denoted by $\mathrm{dom}(D)$,  is the family of square integrable random variables $F=\mathfrak{f}(\Pi)$ such that 
\[
\int_A\mathbb E |D_a F|^2\alpha(\mathrm{d}a)<\infty.
\]

Setting $F:=I(f)$,  we have that $\E F = 0$ with $\mathbb{V}\mathrm{ar}(F)=1$ (as follows by applying the isometry formula for Poisson chaoses) and that $F\in\mathrm{dom}(D)$. 
%$D_a F=f(a)$ and $D_{a_1,a_2}F=0$, $a_1,a_2\in A$.  
%Note also that $F$ is centered with $\mathbb{V}\mathrm{ar}(F)=1$ (as follows by applying the isometry formula for Poisson chaoses) and that
%$F\in\mathrm{dom}(D)$.
 Using Theorem 1.1 in \cite{LPS} we have that 
$$d_W(F,G) \leq \gamma_1 + \gamma_2 + \gamma_3,$$
where 
%$G$ is a standard Gaussian random variable and 
\begin{equation*}
    \begin{cases}
        \gamma_1 &= 2 \left( \int_{A^3}  \left(\E [(D_{a_1}F)^2 (D_{a_2}F)^2]  \right)^{1/2} \left( \E [(D^2_{a_1,a_3}F)^2 (D^2_{a_2,a_3}F)^2]\right)^{1/2} \alpha^{3}(\mathrm {d}(a_1,a_2,a_3))\right)^{1/2} \\
        \gamma_2 &=\left(\int_{A^3} \E [(D^2_{a_1,a_3}F)^2 (D^2_{a_2,a_3}F)^2]\alpha^{3}(\mathrm {d}(a_1,a_2,a_3)) \right)^{1/2}\\
        \gamma_3 &= \int_A \E|D_{a}F|^3 \alpha(\mathrm d a)
    \end{cases}
    .
\end{equation*}
Since for $F=I(f)$ we have that $D_a F=f(a)$ and $D_{a_1,a_2}F=0$, $a_1,a_2\in A$, the terms $\gamma_1$ and $\gamma_2$ vanish and $\gamma_3= \int_A |f(a)|^3 \alpha( \mathrm d a)$.  Therefore,  we obtain the bound on the  Wasserstein distance between $I(f)$ and $G$.
%The bound on the Wasserstein distance is an immediate consequence of Theorem 1.1 in \cite{LPS}. 
Similarly,  Theorem 1.2 in \cite{LPS} gives
%a first upper bound on the Kolmogorov distance between $F$ and $G$ 
$$d_K(F,G) \leq \left( 1 + \frac{1}{2} \left(\E F^4\right)^{1/4} \right) \int_A |f(a)| ^{3} \alpha(\mathrm d a) + \left(\int_A f(a)^4 \alpha(\mathrm d a) \right)^{1/2}.$$
 Using Lemma 4.2 in \cite{LPS} we have
 \begin{align*}
     \E F^4 & \leq \max \left \{256 \left(\int_A f(a)^2 \alpha(\mathrm d a)\right)^{2},  4 \int_A f(a)^4 \alpha(\mathrm d a )+2 \right \}\\
     &= \max \left \{4^4,  4 \int_A f(a)^4 \alpha(\mathrm d a )+2 \right \} \quad \text{(because $\|f\|_{L^2(\alpha)}=1$),}
 \end{align*}
which yields the upper bound on the Kolmogorov distance. 
\\
\noindent$\square$

\subsection{Moderate deviations,  Bernstein-type concentration inequalities and Normal approximation bounds with Cram\'er correction term for first chaoses on the Poisson space}

We start with a definition. 

\begin{Definition}\label{def:27102023}
Let $\{Y_\ell\}_{\ell\in\mathbb N}$ be a sequence of real-valued random variables,  $\gamma\geq 0$ a non-negative constant and 
$\{\Delta_\ell\}_{\ell\in\mathbb N}$ a positive numerical sequence.  We say that:\\
\noindent $(1)$ The sequence $\{Y_\ell\}_{\ell\in\mathbb N}$ satisfies a moderate deviation principle with parameters 
$\gamma$ and $\{\Delta_\ell\}_{\ell\in\mathbb N}$
($\bold{MDP}(\gamma,\{\Delta_\ell\}_{\ell\in\mathbb N})$ for short) if,  for any sequence of positive numbers $\{a_\ell\}_{\ell\in\mathbb N}$ such that $\lim_{\ell\to\infty}a_\ell=+\infty$ and $\lim_{\ell\to\infty}\frac{a_{\ell}}{\Delta_\ell^{1/(1+2\gamma)}}=0$,  the sequence $\{Y_\ell\}_{\ell\in\mathbb N}$ satisfies a large deviation principle with speed $a_\ell^2$ and rate function $J(x):=x^2/2$, i.e., for any Borel set $B\subset\R$,
\[
-\inf_{x\in\overset{\circ}B}J(x)\leq\liminf_{\ell\to\infty}a_\ell^{-2}\log\mathbb{P}\left(Y_\ell/a_\ell\in B\right)\leq
\limsup_{\ell\to\infty}a_\ell^{-2}\log\mathbb{P}\left(Y_\ell/a_\ell\in B\right)\leq-\inf_{x\in\overline B}J(x),
\]
where $\overset{\circ}B$ denotes the interior of $B$ and $\overline B$ denotes the closure of $B$.
\\
\noindent $(2)$ The sequence $\{Y_\ell\}_{\ell\in\mathbb N}$ satisfies a Bernstein-type concentration inequality
with parameters $\gamma$ and $\{\Delta_\ell\}_{\ell\in\mathbb N}$
($\bold{BCI}(\gamma,\{\Delta_\ell\}_{\ell\in\mathbb N})$ for short) if,  for all $\ell\in\mathbb N$ and $x\geq 0$,  we have 
\[
\mathbb{P}(|Y_\ell |\geq x)\leq 2\exp\left(-\frac{1}{4}\min\left\{\frac{x^2}{2^{1+\gamma}},(x\Delta_\ell)^{1/(1+\gamma)}\right\}\right).
\]
\noindent $(3)$ The sequence $\{Y_\ell\}_{\ell\in\mathbb N}$ satisfies a Normal approximation bound with Cram\'er correction term
with parameters $\gamma$ and $\{\Delta_\ell\}_{\ell\in\mathbb N}$
($\bold{NACC}(\gamma,\{\Delta_\ell\}_{\ell\in\mathbb N})$ for short),  if there exist positive constants $c_0,c_1,c_2>0$ only depending on $\gamma$ such that for all $\ell\in\mathbb N$ and $x\in [0,c_0\Delta_\ell^{1/(1+2\gamma)}]$,
\[
\mathbb{P}(Y_\ell\geq x)=\mathrm{e}^{L_{\ell,x}^+}[1-\mathbb P(G\leq x)]\left(1+c_1\theta_{\ell,x}^+\frac{1+x}{\Delta_\ell^{1/(1+2\gamma)}}\right)
\]
and
\[
\mathbb{P}(Y_\ell\leq -x)=\mathrm{e}^{L_{\ell,x}^-}[1-\mathbb P(G\leq -x)]\left(1+c_1\theta_{\ell,x}^-\frac{1+x}{\Delta_\ell^{1/(1+2\gamma)}}\right),
\]
where $\theta_{\ell,x}^{\pm}\in [-1,1]$ and $L_{\ell,x}^{\pm}\in\left(-c_2\frac{x^3}{\Delta_\ell^{1/(1+2\gamma)}},c_2\frac{x^3}{\Delta_\ell^{1/(1+2\gamma)}}\right)$.
\end{Definition}

As a preliminary result, we provide moderate deviations,  Bernstein-type concentration inequalities and Normal approximation bounds with Cram\'er correction term for sequences of first chaoses on the Poisson space.

\begin{Theorem}\label{thm:SchulteThale}
Assume: $(i)$ $f_\ell\in L^m(\alpha_\ell)$ for any $m\geq 2$ and any $\ell\in\mathbb N$ with $\|f_\ell\|_{L^2(\alpha_\ell)}=1$ for any $\ell\in\mathbb N$; $(ii)$ there exist a constant $\gamma\geq 0$ and a positive numerical sequence
$\{\Delta_\ell\}_{\ell\in\mathbb N}$ such that
\[
\Big|\int_A f_\ell(a)^m\alpha_\ell(\mathrm{d}a)\Big|
\leq\frac{(m!)^{1+\gamma}}{\Delta_\ell^{m-2}},\quad\text{for all $m\geq 3$ and $\ell\in\mathbb N$.}
\]

Then the sequence $\{I^{(\ell)}(f_\ell)\}_{\ell\in\mathbb N}$ satisfies a $\bold{MDP}(\gamma,\{\Delta_\ell\}_{\ell\in\mathbb N})$,  a
$\bold{BCI}(\gamma,\{\Delta_\ell\}_{\ell\in\mathbb N})$ and a\\ $\bold{NACC}(\gamma,\{\Delta_\ell\}_{\ell\in\mathbb N})$.

\end{Theorem}

\noindent{\it Proof.}\\
We recall that for real-valued random variables $X_1,\ldots,X_m$,  $m\in\mathbb N$,  the joint cumulant of $X_1,\ldots,X_m$ is defined as
\[
\mathrm{cum}(X_1,\ldots,X_m):=(-\bold{i})^m\frac{\partial^m}{\partial t_1\ldots\partial t_m}\log\varphi_{X_1,\ldots,X_m}(t_1,\ldots, t_m)\Big|_{t_1=\ldots=t_m=0},
\]
where $\bold i$ is the imaginary unit and $\varphi_{X_1,\ldots,X_m}$ is the join characteristic function of $(X_1,\ldots,X_m)$. For a real-valued random variable $X$ and $m\in\mathbb N$ we shall write $\mathrm{cum}_m(X):=\mathrm{cum}(X,\ldots,X)$ for the $m$-th cumulant of $X$.

For an arbitrarily fixed $\ell\in\mathbb N$,  set $X_\ell:=I^{(\ell)}(f_\ell)$. Clearly,
$\mathbb E X_\ell=0$ and $\mathbb E X_\ell^2=1$ (which is a consequence of the isometry formula for Poisson chaoses).  Then the claim follows by 
the theory developed in \cite{SS} (see e.g.  Proposition 2.1 in \cite{ST}; see also \cite{DE} and \cite{DJS}) if we prove that
\begin{equation}\label{eq:22092023}
\mathbb E |X_\ell|^m<\infty, \quad\text{for any $\ell\in\mathbb N$ and $m\geq 3$}
\end{equation}
and

\begin{equation}\label{eq:20092023primo}
\mathrm{cum}_m(X_\ell)=\int_A f_\ell(a)^m\alpha_\ell(\mathrm{d}a),
\quad\text{for any $\ell\in\mathbb N$ and $m\geq 3$.}
\end{equation}
To this aim we are going to apply Theorem 3.6 in \cite{ST}.
A partition $\sigma$ of $\{1,\ldots,m\}$, $m\geq 3$,  is a collection $\{B_1,\ldots,B_k\}$ of $1\leq k\leq m$ pairwise disjoint non-empty sets, called blocks,  such that $B_1\cup\ldots\cup B_k=\{1,\ldots,m\}$. The number $k$ of blocks of a partition $\sigma$ is denoted by $|\sigma|$.  Let $J_j:=\{j\}$,
$j\in\{1,\ldots,m\}$.  Letting $\Pi(\bold{1}_m)$,  $\bold{1}_m:=(1,\ldots,1)\in\R^m$,  denote the set of all partitions $\sigma$ of $\{1,\ldots,m\}$ whose blocks $B$ are such that $\mathrm{Card}(B\cap J_{j})\leq 1$ for every $j\in\{1,\ldots,m\}$,  we clearly have that $\Pi(\bold{1}_m)$ is the set of all partitions of $\{1,\ldots,m\}$. Letting $\widetilde{\Pi}(\bold{1}_m)$
denote the set of all partitions $\sigma\in\Pi(\bold{1}_m)$ with $|\sigma|=1$, we clearly have $\widetilde\Pi(\bold{1}_m)=\{\{1,\ldots,m\}\}$,  $m\geq 3$.
Letting $\widetilde{\Pi}_{\geq 2}(\bold{1}_m)$
denote the set of all partitions $\sigma\in\widetilde\Pi(\bold{1}_m)$ whose blocks have cardinality bigger than or equal to $2$,  since $m\geq 3$,
we clearly have $\widetilde\Pi_{\geq 2}(\bold{1}_m)=\widetilde\Pi(\bold{1}_m)=\{\{1,\ldots,m\}\}$.  We denote by $\Pi_{\geq 2}(\bold{1}_m)$,  $m\geq 2$,  the family of all partitions $\sigma\in\Pi(\bold{1}_m)$, i.e.,  the family of all partitions $\sigma$ of $\{1,\ldots,m\}$,  whose blocks have cardinality bigger than or equal to $2$.

For a function $g:A\to\mathbb R$, set
\[
(\otimes_{j=1}^{m}g)(x_1,\ldots,x_m):=g(x_1)\ldots g(x_m).
\]
For $\sigma\in\Pi(\bold{1}_m)$,  define the function $(\otimes_{j=1}^{m}g)_\sigma:A^{|\sigma|}\to\R$ by replacing in $(\otimes_{j=1}^{m}g)(x_1,\ldots,x_m)$
all the variables whose indexes belong to the same block of $\sigma$ by a new common variable.  Note that for $\sigma\in\Pi(\bold{1}_m)$,  $(\otimes_{j=1}^{m}g)_\sigma:A^{|\sigma|}\to\R$
can be represented as

\[
(\otimes_{j=1}^{m}g)_\sigma(a_1,\ldots,a_{|\sigma|})=\prod_{i=1}^{|\sigma|}g(a_i)^{|B_i|},\quad a_1,\ldots,a_{|\sigma|}\in A
\]
where $B_1,\ldots,B_{|\sigma|}$ are the blocks of $\sigma$ and $|B_i|$, $i=1,\ldots,|\sigma|$, is the cardinality of the block $B_i$.
In particular, for $\sigma\in\widetilde\Pi_{\geq 2}(\bold{1}_m)$,  
\[
(\otimes_{j=1}^{m}g)_\sigma(a):=g(a)^m,\quad a\in A
\]
and,  for $\sigma\in\Pi _{\geq 2}(\bold{1}_m)$,  
\[
(\otimes_{j=1}^{m}g)_\sigma(a_1,\ldots,a_{|\sigma|}):=\prod_{i=1}^{|\sigma|}g(a_i)^{|B_i|},\quad a_1,\ldots,a_{|\sigma|}\in A
\]
where $B_1,\ldots,B_{|\sigma|}$ are the blocks of $\sigma$ and $|B_i|\geq 2$ for any $i=1,\ldots,|\sigma|$. Therefore the hypothesis $(i)$
implies the assumptions $(3.4)$ and $(3.5)$ of Theorem 3.6 in \cite{ST},  and so
\eqref{eq:22092023} and \eqref{eq:20092023primo} hold.\\
\noindent$\square$

\section{Application to Poisson shot noise random variables}\label{sec:APPPSN}

\subsection{Gaussian approximation}

In this section, we apply Theorems \ref{thm:LPS1} to the standardized random variables $T(C)$,  $C\in\mathcal{B}(\mathbb R^d)$,  defined by \eqref{eq:11102023terzo}.
Note that, for $C\in\mathcal{B}(\R^d)$,
\[
\mathbb E S(C)=\int_{\mathbb R^d} \lambda(x) \mathbb E H(C-x,Z_1) \mathrm d x,
\]
and that,  by the isometry formula for Poisson chaoses, 
if $\int_{\mathbb R^d}\lambda(x)\mathbb{E}H(C-x,Z_1)^2\mathrm{d}x<\infty$, then 
\[
\mathbb{V}\mathrm{ar}(S(C))=\int_{\mathbb R^d}\lambda(x)\mathbb{E}H(C-x,Z_1)^2\mathrm{d}x.
\]
Therefore the random variable
$T(C)$ is well defined and finite for any $C \in \mathcal{B}(\mathbb R^d)$ such that 
\begin{equation}\label{eq:intzeroinfty}
0 <\int_{\mathbb R^d}\lambda(x)\mathbb{E}H(C-x,Z_1)^2\mathrm{d}x<\infty.
\end{equation}

The following theorem holds.

\begin{Theorem}\label{thm:GaussPoissoncluster}
Let $C \in \mathcal{B}(\mathbb R^d)$ be such that \eqref{eq:intzeroinfty} holds.

Then
\[
d_W(T(C),G)\leq
\frac{\int_{\R^d}\lambda(x)\mathbb{E}|H(C-x,Z_1)|^3\mathrm{d}x}{\left(\int_{\R^d}\lambda(x)\mathbb{E}H(C-x,Z_1)^2\mathrm{d}x\right)^{3/2}}
\]
and
\begin{align}
&d_K(T(C),G)\nonumber\\
&\qquad
\leq
\left[1+\frac{1}{2}\max\Biggl\{4,\left[4\frac{\int_{\R^d}\mathbb \lambda(x)\mathbb{E}H(C-x,Z_1)^4\mathrm{d}x}{\left(\int_{\R^d} \lambda(x) \mathbb{E}H(C-x,Z_1)^2\mathrm{d}x\right)^2}+2\right]^{1/4}\Bigg\}\right]\frac{\int_{\R^d}\lambda(x)\mathbb{E}|H(C-x,Z_1)|^3\mathrm{d}x}{\left(\int_{\R^d}\lambda(x)\mathbb{E}H(C-x,Z_1)^2\mathrm{d}x\right)^{3/2}}\nonumber\\
&\qquad\qquad\qquad
+\left(\frac{\int_{\R^d}\mathbb \lambda(x)\mathbb{E}H(C-x,Z_1)^4\mathrm{d}x}{\left(\int_{\R^d}  \lambda(x)\mathbb{E}H(C-x,Z_1)^2\mathrm{d}x\right)^2}\right)^{1/2}.
\nonumber
\end{align}
\end{Theorem}

\begin{Remark}\label{re:particularcases}
As particular cases of Theorem \ref{thm:GaussPoissoncluster} we have the following:\\
\noindent $(i)$ If $H(C-x,z):=v(z)(C-x)$,  with $C\in\mathcal{B}(\mathbb R^d)$,  $x\in\mathbb{R}^d$, $z\in\bold{Z}:=\bold{N}_{\mathbb{R}^d\times\mathbb{R}^p}$
and $\upsilon(z)(C)$ is defined by \eqref{eq:15122023pom5},  
then $S(C)=V(C)$ where the random variable $V(C)$ is defined by \eqref{eq:rapprS}.  So Theorem \ref{thm:GaussPoissoncluster} provides Gaussian approximation bounds for the random variable $W(C)$ defined by \eqref{eq:11102023terzoV}.  An interesting particular case is investigated in Section \ref{sec:1512partic1}.\\
\noindent $(ii)$ If $H:\mathcal{B}(\R^2)\times (0,\infty)\to (0,\infty)$ is a mapping such that its restriction on $\R^2\times (0,\infty)$,
coincides with $\widetilde{H}(x,z):=zA(-x)$,  $x\in\mathbb R^2$ and $z\in\bold{Z}:=(0,\infty)$,  then
then $S(\{\bold 0\})=I(\{\bold 0\})$ where the random variable $I(\{\bold 0\})$ is defined in Section \ref{sec:07122023}.  So Theorem \ref{thm:GaussPoissoncluster} provides Gaussian approximation bounds for the random variable $L(\{\bold 0\})$ defined by \eqref{eq:11102023terzoI}.  An interesting particular case is investigated in Section
\ref{sec:1512partic3}.

\end{Remark}

\noindent{\it Proof\,\,of\,\,Theorem\,\,\ref{thm:GaussPoissoncluster}.}\\
By \eqref{eq:PSN} we have
\begin{equation*}
T(C)=\frac{1}{\sqrt{\int_{\mathbb R^d}\lambda(x)\mathbb{E}H(C-x,Z_1)^2\mathrm{d}x}}\int_{\mathbb R^d\times\bold{M}}H(C-x,z)(\mathcal{P}(\mathrm{d}x,\mathrm{d}z)-\lambda(x)\mathrm{d}x\mathbb{Q}(\mathrm{d}z)),\quad C\in\mathcal{B}(\R^d).
\end{equation*}
Therefore $T(C)$ belongs to the first chaos of $\mathcal P$ with kernel
\[
t(x,z):=\frac{H(C-x,z)}{\sqrt{\int_{\mathbb R^d}\lambda(x)\mathbb{E}H(C-x,Z_1)^2\mathrm{d}x}}.
\] 
The claim follows applying Theorem \ref{thm:LPS1}.\\

\noindent$\square$

\subsection{Moderate deviations,  Bernstein-type concentration inequalities and Normal approximation bounds with Cram\'er correction term}

In this section, we apply Theorem \ref{thm:SchulteThale} to the sequence of standardized random variables $\{T_\ell(C_\ell)\}_{\ell\geq 1}$,  $\{C_\ell\}_{\ell\geq 1}\subset\mathcal{B}(\mathbb R^d)$,  defined by \eqref{eq:PSNsequence}.

The following theorem holds.

\begin{Theorem}\label{thm:ModeratePoissoncluster}
Let $\{C_\ell\}_{\ell \in \mathbb N}\subset\mathcal{B}(\mathbb R^d)$ be a sequence of Borel sets such that 
\begin{equation}\label{eq:29102023I}
0<\int_{\R^d}\lambda_\ell(x)\mathbb{E}H(C_\ell-x,Z_1^{(\ell)})^2\mathrm{d}x<\infty,\quad \ell\geq 1
\end{equation}
and assume that there exist a non-negative constant $\gamma\geq 0$ and a positive numerical sequence $\{\Delta_\ell\}_{\ell \in \mathbb N}$ such that 
\begin{equation}\label{eq:08122023primo}
\frac{\int_{\R^d} \lambda_\ell (x) \E |H(C_\ell-x,Z_1^{(\ell)})|^m \mathrm d x}{\left( \int_{\R^d} \lambda_\ell (x) \E H(C_\ell-x,Z_1^{(\ell)})^2 \mathrm d x\right)^{\frac{m}{2}}} \leq \frac{(m!)^{1+\gamma}}{\Delta_\ell ^{m-2}},\quad \text{for all } m\geq 3 \text{ and }\ell\in\mathbb N.
\end{equation}
Then the sequence $\{T_\ell(C_\ell)\}_{\ell\geq 1}$ satisfies a $\bold{MDP}(\gamma,\{\Delta_\ell\}_{\ell\in\mathbb N})$,  a $\bold{BCI}(\gamma,\{\Delta_\ell\}_{\ell\in\mathbb N})$ and a\\ $\bold{NACC}(\gamma,\{\Delta_\ell\}_{\ell\in\mathbb N})$.
\end{Theorem}

\begin{Remark}\label{re:particularcasesbis}
As particular cases of Theorem \ref{thm:ModeratePoissoncluster} we have the following:\\ 
\noindent $(i)$ If $H(C_\ell-x,z):=v(z)(C_\ell-x)$,  with $C_\ell\in\mathcal{B}(\mathbb R^d)$,  $x\in\mathbb{R}^d$, $z\in\bold{Z}:=\bold{N}_{\mathbb{R}^d\times\mathbb{R}}$ and $\upsilon(z)(\cdot)$ is defined by \eqref{eq:15122023pom5},  
then $S_\ell(C_\ell)=V_\ell(C_\ell)$ where the random variable $V_\ell(C_\ell)$ is defined by \eqref{eq:rapprS} with $\lambda_\ell$,  $\mathbb Q_\ell$ and $C_\ell$
in place of $\lambda$,  $\mathbb Q$ and $C$,  respectively.  So Theorem \ref{thm:ModeratePoissoncluster} provides moderate deviations,  Bernstein-type concentration inequalities and Normal approximation bounds with Cram\'er correction term for the sequence $\{W_\ell(C_\ell)\}_{\ell\geq 1}$,  where $W_\ell(C_\ell)$ is defined by \eqref{eq:11102023terzoV}, with obvious modifications.  An interesting particular case is investigated in Section \ref{sec:1512partic2}.\\
\noindent $(ii)$ If $H:\mathcal{B}(\R^2)\times (0,\infty)\to (0,\infty)$ is a mapping such that its restriction on $\R^2\times (0,\infty)$,
coincides with $\widetilde{H}(x,z):=zA(-x)$,  $x\in\mathbb R^2$ and $z\in\bold{Z}:=(0,\infty)$,  then
$S_\ell(\{\bold 0\})=I_\ell(\{\bold 0\})$ where the random variable $I_\ell(\{\bold 0\})$ is defined as in Section \ref{sec:07122023},  with $\lambda_\ell$ and $\mathbb Q_\ell$
in place of $\lambda$ and $\mathbb Q$,  respectively.  So Theorem \ref{thm:ModeratePoissoncluster} provides moderate deviations,  Bernstein-type concentration inequalities and Normal approximation bounds with Cram\'er correction term for the sequence $\{L_\ell(\{\bold 0\})\}_{\ell\geq 1}$,  where $L_\ell(\{\bold 0\})$ is defined by \eqref{eq:11102023terzoI}, with obvious modifications.  An interesting particular case is investigated in Section \ref{sec:1512partic4}.
\end{Remark}

\noindent{\it Proof\,\,of\,\,Theorem\,\,\ref{thm:ModeratePoissoncluster}.}\\ 
Similarly to the proof of Theorem \ref{thm:GaussPoissoncluster},
by \eqref{eq:PSNsequence} we have
\begin{equation*}
T_\ell(C_\ell)=\frac{1}{\sqrt{\int_{\mathbb R^d}\lambda_\ell(x)\mathbb{E}H(C_\ell-x,Z_1^{(\ell)})^2\mathrm{d}x}}\int_{\mathbb R^d\times\bold{Z}}H(C_\ell-x,z)(\mathcal{P}_\ell(\mathrm{d}x,\mathrm{d}z)-\lambda_\ell(x)\mathrm{d}x\mathbb{Q}_\ell(\mathrm{d}z)),\quad C_\ell\in\mathcal{B}(\R^d).
\end{equation*}
Therefore $T_\ell(C_\ell)$ belongs to the first chaos of $\mathcal P_\ell$ with kernel
\[
t_\ell(x,z):=\frac{H(C_\ell-x,z)}{\sqrt{\int_{\mathbb R^d}\lambda_\ell(x)\mathbb{E}H(C_\ell-x,Z_1^{(\ell)})^2\mathrm{d}x}}.
\] 
The claim follows by Theorem \ref{thm:SchulteThale}. \\

\noindent$\square$

\section{Application to a class of compound Poisson cluster point processes}\label{sec:CPCP}

\subsection{Gaussian approximation}\label{sec:1512partic1}

In this section we apply Theorem \ref{thm:GaussPoissoncluster} to the class of standardized compound Poisson cluster point processes $\{W(C)\}_{C\in\mathcal{B}(\R^d)}$ 
with $\{X_n\}_{n\geq 1}$ having a piecewise constant intensity function.  In such a case we have more explicit upper bounds on the Wasserstein and the Kolmogorov distances,  which pave the way to explicit bounds 
for some classes of generalized compound Hawkes processes (see Corollaries \ref{cor:11102023quarto} and \ref{cor:11102023quartoBIS}).  We recall that $Z$ denotes the random total number of points of the progeny process, i.e.  $Z=Z_1(\R^d,\R)$,  and that $M$ is a generic random variable with the same distribution as the independent and identically distributed marks.

Hereafter,  we denote by $\mathrm{Leb}(\cdot)$ the Lebesgue measure on $\mathbb R^d$.

\begin{Corollary}\label{cor:2}
Let $(B,C)\in\mathcal{B}(\mathbb R^d)^2$ be such that $0<\mathrm{Leb}(B\cap C)<+\infty$.  If $\lambda(x)=\lambda \boldsymbol 1 _{B}(x)$ 
for any $x\in\R^d$ and some positive constant $\lambda>0$,  $\mathbb E Z^2<\infty$ and $\E M^2\in (0,\infty)$,  then

\begin{equation}\label{eq:11102023primo}
d_W(W(C),G)\leq \frac{\E |M|^3\mathbb E Z^3}{(\E M^2)^{3/2}\sqrt {\lambda \mathrm{Leb}(B\cap C)}}
\end{equation}
and 
\begin{align}
d_K(W(C),G)
&\leq
\left[1+\frac{1}{2}\max\Biggl\{4,\left[4\frac{\E M^4\mathbb{E} Z^4}{\lambda\mathrm{Leb}(B\cap C)(\E M^2)^2}+2\right]^{1/4}\Bigg\}\right]\frac{\E |M|^3\mathbb E Z^3}{(\E M^2)^{3/2}\sqrt {\lambda \mathrm{Leb}(B\cap C)}}\nonumber\\
&\qquad\qquad\qquad\qquad
+\left(\frac{\E M^4\mathbb{E} Z^4}{\lambda \mathrm{Leb}(B\cap C)(\E M^2)^2}\right)^{1/2}.
\label{eq:11102023secondo}
\end{align}
\end{Corollary}

 We point out  that many articles in the literature (\textit{e.g.} \cite{Bacry,HHKR,KPR}) consider Hawkes processes with an empty history, that is with no points in $(-\infty,0]$, which corresponds to the piecewise constant intensity function $\lambda(x)=\boldsymbol 1_{[0,+\infty)}(x)$ in Corollary \ref{cor:2}.

The proof of Corollary \ref{cor:2} exploits the following lemma,  which is proved
in Section \ref{sec:lemmas}.

\begin{Lemma}\label{lemma:Z}
For any $(B,C)\in\mathcal{B}(\mathbb R^d)^2$ and $m\in\mathbb N$,  we have
\begin{equation}\label{eq:23052023settimo}
\int_{B} \E |\upsilon(Z_1) (C-x)| ^m \mathrm d x \leq \mathrm{Leb}(B\cap C)\E Z^m\E |M|^m .
\end{equation}
\end{Lemma}

\noindent{\it Proof\,\,of\,\,Corollary\,\,\ref{cor:2}.}\\
In order to apply Theorem \ref{thm:GaussPoissoncluster} we need to verify \eqref{eq:intzeroinfty} with 
$H(C-x,z):=\upsilon(z)(C-x)$,  $C\in\mathcal{B}(\mathbb R^d)$,  $x\in\mathbb R^d$, $z\in\bold{Z}:=\bold{N}_{\mathbb R^d\times\R}$,
and $\lambda(\cdot)\equiv\lambda \boldsymbol 1 _B(\cdot)$.  For the lower bound we note that
    \begin{align*}
        \int_{B} \mathbb E \upsilon(Z_1) (C-x)^2 \mathrm d x=\int_{B}\mathbb E\left(\sum _{k=0}^{Z_1(C-x,\R)-1}M_{1,k}\right)^2\mathrm d x
         \geq \int_{B\cap C}\mathbb E\left(\sum _{k=0}^{Z_1(C-x,\R)-1}M_{1,k}\right)^2\mathrm d x.
    \end{align*}
    Expanding the square of the sum and using the independence we have that 
    $$\mathbb E\left(\sum _{k=1}^{Z_1(C-x,\R)}M_{1,k}\right)^2=\E[Z_1(C-x,\R)]\E M^2+\E[Z_1(C-x,\R)(Z_1(C-x,\R)-1)](\E M)^2\nonumber.$$
    
    For $x\in B\cap C$,  we have $x\in C$,  and so the set $C-x$ contains the origin.  Since $Z_1(\{\bold 0\},\R)=1$,  we then have
    $$\mathbb E\left(\sum _{k=1}^{Z_1(C-x,\R)}M_{1,k}\right)^2 \geq \E M^2.$$
    Therefore,
    \begin{align}
        \int_{B} \mathbb E \upsilon(Z_1) (C-x)^2 \mathrm d x & \geq 
\int_{B\cap C} \mathbb E \upsilon(Z_1) (C-x)^2 \mathrm d x\geq 
        \mathrm{Leb}(B\cap C)\mathbb E M^2>0.\label{eq:23052023ottavo}
    \end{align}

As far as the upper bound is concerned,  we note that by the bound on the Wasserstein distance in Theorem \ref{thm:GaussPoissoncluster} and the inequalities \eqref{eq:23052023ottavo} and \eqref{eq:23052023settimo} immediately follows
    \begin{align*}
        d_W(F(C),G)&\leq
\frac{\lambda\int_{B}\mathbb{E}|\upsilon(Z_1)(C-x)|^3\mathrm{d}x}{\left(\int_{B}\lambda\mathbb{E}\upsilon(Z_1)(C-x)^2\mathrm{d}x\right)^{3/2}}
\leq \frac{\lambda \mathrm{Leb}(B\cap C)\mathbb E Z^3 \E |M|^3}{\left(\lambda\E M^2\mathrm{Leb}(B\cap C) \right)^{3/2}}
=\frac{\E |M|^3\mathbb E Z^3}{(\E M^2)^{3/2}\sqrt{\lambda \mathrm{Leb} (B\cap C)}}.
    \end{align*}
Similarly,  the upper bound on the Kolmogorov distance follows by the upper bound on the Kolmogorov distance in Theorem \ref{thm:GaussPoissoncluster}, and again
the inequalities \eqref{eq:23052023ottavo} and \eqref{eq:23052023settimo}.
\\
\noindent$\square$

\subsection{Moderate deviations,  Bernstein-type concentration inequalities and Normal approximation bounds with Cram\'er correction term}
\label{sec:1512partic2}

In this section we apply Theorem \ref{thm:ModeratePoissoncluster} to the sequence $\{W_\ell(C_\ell)\}_{\ell\geq 1}$,  when the Poisson processes $\{X_n^{(\ell)}\}$,  $\ell\geq 1$, have a piecewise deterministic intensity function and $\mathbb Q_\ell\equiv\mathbb Q$,  for every $\ell\geq 1$. 
  
In such a case the assumption \eqref{eq:08122023primo} simplifies a lot.  Moreover,  the next corollary paves the way to the application to some classes of generalized compound Hawkes processes (see Corollaries \ref{cor:11102023quinto} and \ref{cor:11102023nono}).

\begin{Corollary}\label{cor:ST}
Let $\{(B_\ell,C_\ell)\}_{\ell\in\mathbb N}\subset\mathcal{B}(\mathbb R^d)^2$ 
be such that $0<\mathrm{Leb}(B_\ell \cap C_\ell)<+\infty$,  $\ell\in\mathbb N$.  Assume that
$\lambda_\ell(x)=\lambda_\ell \boldsymbol 1_{B_\ell}(x)$,  $x\in\mathbb R^d$,  for  
positive constants $\lambda_\ell>0$,  $\ell\in\mathbb N$,
$\mathbb Q_\ell\equiv \mathbb Q$,  $\E M^2>0$ and
\begin{equation}\label{hyp:29092023}
\frac{\E |M|^m\E Z^m}{(\E M^2)^{m/2}\sqrt{\lambda_\ell \mathrm{Leb}(B_\ell \cap C_\ell)}^{m-2}}\leq \frac{(m!)^{1+\gamma}}{\Delta_\ell ^{m-2}},\quad\text{for all $m\geq 3$ and $\ell\in\mathbb N$}
\end{equation}
for some $\gamma\geq 0$ and a positive numerical sequence $\{\Delta_\ell\}_{\ell \in \mathbb N}$.  Then the sequence $\{W_\ell(C_\ell)\}_{\ell\geq 1}$ satisfies a
$\bold{MDP}(\gamma,\{\Delta_\ell\}_{\ell\in\mathbb N})$,  a $\bold{BCI}(\gamma,\{\Delta_\ell\}_{\ell\in\mathbb N})$ and a $\bold{NACC}(\gamma,\{\Delta_\ell\}_{\ell\in\mathbb N})$.

\end{Corollary}

\noindent{\it Proof\,\,of\,\,Corollary\,\,\ref{cor:ST}.}\\ 
By \eqref{eq:23052023ottavo}, 
\eqref{eq:23052023settimo},  the choice of the Borel sets $B_\ell$ and $C_\ell$,  $\ell\in\mathbb N$,  the assumption \eqref{hyp:29092023} and the fact that $\E M^2>0$, 
we have
\begin{equation*}
0<\mathrm{Leb}(B_\ell \cap C_\ell)\E M^2\leq\int_{B_\ell} \mathbb E\upsilon(Z_1)(C_\ell-x)^2 \mathrm d x
\end{equation*}
and
\begin{equation*}
\int_{B_\ell} \mathbb E|\upsilon(Z_1)(C_\ell-x)|^m \mathrm d x \leq\mathrm{Leb}(B_\ell \cap C_\ell)\E |M|^m\mathbb E Z^m<\infty,\quad\text{for every $m\geq 3$.}
\end{equation*}
Therefore condition \eqref{eq:29102023I} holds and using again the assumption \eqref{hyp:29092023},

we have
\begin{align*}
\frac{\lambda_\ell\int_{B_\ell}  \E |\upsilon(Z_1)(C_\ell-x)|^m \mathrm d x}{\left( \lambda_\ell\int_{B_\ell}  \E (\upsilon(Z_1)(C_\ell-x))^2 \mathrm d x\right)^{\frac{m}{2}}}&\leq \frac{\lambda_\ell \mathrm{Leb}(B_\ell \cap C_\ell)\E |M|^m\mathbb E Z^m}{(\lambda_\ell \mathrm {Leb}(B_\ell \cap C_\ell)\E M^2)^{\frac{m}{2}} }\\ &= \frac{\E |M|^m\mathbb E Z^m}{(\E M^2)^{m/2}\sqrt{\lambda_\ell \mathrm{Leb}(B_\ell \cap C_\ell)}^{m-2}}\\ &\leq \frac{(m!)^{1+\gamma}}{\Delta_\ell ^{m-2}},\quad\text{for all $m\geq 3$ and $\ell\in\mathbb N$.}
\end{align*}
The claim follows by Theorem \ref{thm:ModeratePoissoncluster}.
\\
\noindent$\square$

\section{Application to generalized compound Hawkes processes}\label{sec:CHP}

We start with a proposition which expresses the moments of the total progeny of a Galton-Watson process with one ancestor in terms of moments of the offspring distribution.

\begin{Proposition}\label{thm:branching}
Suppose that $Z$ is distributed as the total progeny of a Galton-Watson process with one ancestor and let $P$ be a random variable 
distributed according to the offspring law of the Galton-Watson process. Assume that the Galton-Watson process is subcritical, i.e.
\begin{equation}\label{eq:sub}
\mathbb E P\in (0,1). 
\end{equation}
Then, for any $n\in \mathbb N$ such that 
\begin{equation}\label{eq:2023}
\mathbb E [P^n]<+\infty,
\end{equation}
we have that $\mathbb E Z^n<\infty$ and 

\begin{align}
\mathbb{E}[Z^n]&=1+
\sum_{k=1}^{n}k!\binom{n}{k}
\sum_{i=1}^{k}\frac{\mathbb{E}[(P)_i]}{i!}
\sum_{m_1+m_2+\ldots+m_i=k}
\frac{\mathbb{E}[Z^{m_1}]}{m_1!}\ldots\frac{\mathbb{E}[Z^{m_i}]}{m_i!},
\label{eq:recursive}
\end{align}
where $\mathbb E(P)_1=\mathbb E P$,
\[
\mathbb{E}(P)_i=\mathbb E P(P-1)\cdot\ldots\cdot P(P-(i-1)),\quad 2\leq i\leq n
\]
and the third sum in \eqref{eq:recursive} is taken over all the $m_1,\ldots,m_i\in\mathbb N$ such that $m_1+\ldots+m_i=k$.

In particular,
\[
\mathbb E Z^2=\frac{\mathbb{V}\mathrm{ar}(P)+1-\mathbb E P}{(1-\mathbb E P)^3},
\]
\begin{align}
&\mathbb E Z^3=\frac{1}{1-\mathbb E P}\Biggl(1+3\frac{\mathbb E P}{1-\mathbb E P}+3\frac{\mathbb E (P)_2}{(1-\mathbb E P)^2}+\frac{\mathbb E (P)_3+3\mathbb{V}\mathrm{ar}(P)}{(1-\mathbb E P)^3}
+3\frac{\mathbb{V}\mathrm{ar}(P)^2}{(1-\mathbb E P)^4}\Biggr)\nonumber
\end{align}
and
\begin{align*}
    \mathbb E Z^4=&\frac{1}{1-\mathbb E P} \Bigg [1+ 4\frac{\mathbb E P}{1-\mathbb E P} +6 \frac{\mathbb E (P)_2}{(1-\mathbb E P)^2}+4 \frac{\mathbb E (P)_3}{(1-\mathbb E P)^3} + \frac{\mathbb E (P)_4}{(1-\mathbb E P)^4 }\\
    &+ 3 \mathbb E Z^2 \Big( 2 \mathbb E P + 4 \frac{\mathbb E (P)_2}{1-\mathbb E P}+ \mathbb E (P)_2 \mathbb E Z^2  +2 \frac{\mathbb E (P)_3}{(1-\mathbb E P)^2}\Big) +4 \mathbb E Z^3 \frac{\mathbb{V}\mathrm{ar}(P)}{1-\mathbb E P} \Bigg].
\end{align*}

\end{Proposition}

As immediate consequence of this proposition and Corollary \ref{cor:2},  we have that if the point processes $Z_i(\cdot,\R)$ are such that $Z$ is distributed as the total progeny of a Galton-Watson process with one ancestor and offspring law satisfying \eqref{eq:sub},  then $(i)$ If the offspring law satisfies 
\eqref{eq:2023} with $n=3$ and $\E M^2\in (0,\infty)$,  then relation
\eqref{eq:11102023primo} holds and the upper bound on $d_W(W(C),G)$ is explicit and depends only on $\lambda$,  the Lebesgue measure of $B\cap C$,  the first three moments of the offspring law and the second and third moment of $|M|$; $(ii)$  If the offspring law satisfies 
\eqref{eq:2023} with $n=4$ and $\E M^2\in (0,\infty)$,  then relation \eqref{eq:11102023secondo} holds and the upper bound on $d_K(W(C),G)$ is explicit and depends only on $\lambda$,  the Lebesgue measure of $B\cap C$,  the first four moments of the offspring law and the second, the third and the fourth moment of $|M|$.

The cases of the Poisson offspring law (which includes compound Hawkes processes) and of the Binomial offspring law are treated in detail in Sections \ref{subsec:gausspoisson}
and \ref{subsec:gaussbinom},  respectively.

The proof of Theorem \ref{thm:branching} exploits the following lemma. Hereafter, for a sufficiently smooth function $f$, we denote by $f^{(n)}$ its derivative of order $n\in\mathbb N$.

\begin{Lemma}\label{le:FaadiBruno} $($Fa\`a\,\,di\,\,Bruno\,\,formula$)$
For any sufficiently smooth functions $g$ and $h$,
\[
(g\circ h)^{(j)}(x)=j!\sum_{i=1}^{j}\frac{g^{(i)}(h(x))}{i!}\sum_{m_1+m_2+\ldots+m_i=j}\frac{h^{(m_1)}(x)}{m_1!}\ldots\frac{h^{(m_i)}(x)}{m_i!},
\quad j\in\mathbb N,
\]
where the second sum is taken over all the $m_1,\ldots,m_i\in\mathbb N$ such that $m_1+\ldots+m_i=j$. 
\end{Lemma}

\noindent{\it Proof\,\,of\,\,Proposition\,\,\ref{thm:branching}}. \\ 
We divide the proof in three steps.  In the first step we provide a functional equation for $\mathbb E\mathrm{e}^{\theta Z}$, $\theta\in (-\infty,0)$.  In the second step we prove the finiteness of the moments of $Z$ and the formula \eqref{eq:recursive}.  In the third step we compute the second, the third and the fourth moments of $Z$.\\
\noindent{\it Step\,\,1: A functional equation for $\mathbb E \mathrm{e}^{\theta Z}$,  $\theta\in(-\infty,0)$.}\\
We note that $Z$ can be represented as 
\[
Z=\sum_{n\geq 0}K_n,
\]
where $K_0=1$ and $K_n$ is the number of offspring in the $n$th generation of the related Galton-Watson process.
Let $\{Z_j\}_{j\geq 1}$ be independent copies of $Z$.
For any $\theta\in(-\infty,0)$, 
by standard computations we have
\begin{align}
\mathbb{E}[\mathrm{e}^{\theta Z}]
&=\mathrm{e}^{\theta}\sum_{k\geq 0}\mathbb{E}[\mathrm{e}^{\theta\sum_{j=1}^{k}Z_j}\,|\,K_1=k]p_k\nonumber\\
&=\mathrm{e}^{\theta}\sum_{k\geq 0}\mathbb{E}[\mathrm{e}^{\theta Z}]^{k}p_k=\mathrm{e}^{\theta}G_{P}(\mathbb{E}[\mathrm{e}^{\theta Z}])<\infty,\label{eq:eqLap}
\end{align}
where $\{p_k\}_{k\geq 0}$ is the law of $P$ and $G_P$ is the probability generating function of $P$.\\
\noindent{\it Step\,\,2: Proof\,\,of\,\,$\mathbb E Z^n<\infty$,  $n\in\mathbb N$,  and of \eqref{eq:recursive}.}\\
As far as the moments of $Z$ are concerned, 
we start by showing that they coincide with the left derivative of the moment generating function at zero. For any $\theta < 0$, the theorem of differentiation under the expected value yields,  for any non-negative integer $n$,
$$\frac{\mathrm d^n}{\mathrm d \theta ^n}\mathbb E[e^{\theta Z}]= \mathbb E [Z^n e^{\theta Z}]<\infty$$
The family $(Z^n e^{\theta Z})_{\theta <0}$ is nonnegative and increasing in $\theta$, hence, using the Beppo Levi theorem
\begin{align}
    \lim_{\theta \uparrow 0} \frac{\mathrm d^n}{\mathrm d \theta ^n}\mathbb E[e^{\theta Z}]&= \lim_{\theta \uparrow 0} \mathbb E \left [ Z^n e^{\theta Z}\right]\nonumber\\
    &=\mathbb E \left [ \lim_{\theta \uparrow 0} Z^n e^{\theta Z}\right]\nonumber\\
    &=\mathbb E \left [  Z^n \right],\label{eq:30052023primo}
\end{align}
where the equality holds whether the quantities are finite or infinite. 
Next, we combine the Fa\`a di Bruno formula with the elementary relation: 
\begin{equation}
\label{eq:fprod}
\frac{\mathrm{d}^n}{\mathrm{d}x^n}(f(x)g(x))=\sum_{k=0}^{n}\binom{n}{k}f^{(n-k)}(x)g^{(k)}(x),\quad n\in\mathbb N\cup\{0\}
\end{equation}
for sufficiently smooth functions $f$ and $g$. For any $\theta\in (-\infty,0)$, by \eqref{eq:eqLap} and \eqref{eq:fprod}, for any non-negative integer $n$,  we have
\[
\frac{\mathrm{d}^n}{\mathrm{d}\theta^n}\mathbb{E}[\mathrm{e}^{\theta Z}]=\mathrm{e}^{\theta}\sum_{k=0}^{n}\binom{n}{k}
\frac{\mathrm{d}^k}{\mathrm{d}\theta^k}G_{P}(\mathbb{E}[\mathrm{e}^{\theta Z}]).
\]
By the Fa\`a di Bruno formula we have
\[
\frac{\mathrm{d}^k}{\mathrm{d}\theta^k}G_{P}(\mathbb{E}[\mathrm{e}^{\theta Z}])
=k!\sum_{i=1}^{k}\frac{G_{P}^{(i)}(\mathbb{E}[\mathrm{e}^{\theta Z}])}{i!}
\sum_{m_1+m_2+\ldots+m_i=k}\frac{
\frac{\mathrm{d}^{m_1}}{\mathrm{d}\theta^{m_1}}
\mathbb{E}[\mathrm{e}^{\theta Z}]}{m_1!}\ldots\frac{\frac{\mathrm{d}^{m_i}}{\mathrm{d}\theta^{m_i}}
\mathbb{E}[\mathrm{e}^{\theta Z}]}{m_i!},
\quad k\in\mathbb N,
\]
where the sum is taken over all the $m_1,\ldots,m_i\in\mathbb N$ such that $m_1+\ldots+m_i=k$.  Then, for any $\theta\in (-\infty,0)$ and $n\in\mathbb N\cup\{0\}$,
\begin{align}
\frac{\mathrm{d}^n}{\mathrm{d}\theta^n}\mathbb{E}[\mathrm{e}^{\theta Z}]&=\mathrm {e}^\theta G_P(\mathbb E\mathrm{e}^{\theta Z})+\mathrm{e}^{\theta}\sum_{k=1}^{n}\binom{n}{k}
\frac{\mathrm{d}^k}{\mathrm{d}\theta^k}G_{P}(\mathbb{E}[\mathrm{e}^{\theta Z}])\nonumber\\
&=\mathrm {e}^\theta G_P(\mathbb E\mathrm{e}^{\theta Z})+
\mathrm{e}^{\theta}\sum_{k=1}^{n}\binom{n}{k}
k!\sum_{i=1}^{k}\frac{G_{P}^{(i)}(\mathbb{E}[\mathrm{e}^{\theta Z}])}{i!}
\sum_{m_1+m_2+\ldots+m_i=k}\frac{
\frac{\mathrm{d}^{m_1}}{\mathrm{d}\theta^{m_1}}
\mathbb{E}[\mathrm{e}^{\theta Z}]}{m_1!}\ldots\frac{\frac{\mathrm{d}^{m_i}}{\mathrm{d}\theta^{m_i}}
\mathbb{E}[\mathrm{e}^{\theta Z}]}{m_i!}\nonumber\\
&=\mathrm {e}^\theta G_P(\mathbb E\mathrm{e}^{\theta Z})+
\mathrm{e}^{\theta}\sum_{k=1}^{n-1}\binom{n}{k}
k!\sum_{i=1}^{k}\frac{G_{P}^{(i)}(\mathbb{E}[\mathrm{e}^{\theta Z}])}{i!}
\sum_{m_1+m_2+\ldots+m_i=k}\frac{
\frac{\mathrm{d}^{m_1}}{\mathrm{d}\theta^{m_1}}
\mathbb{E}[\mathrm{e}^{\theta Z}]}{m_1!}\ldots\frac{\frac{\mathrm{d}^{m_i}}{\mathrm{d}\theta^{m_i}}
\mathbb{E}[\mathrm{e}^{\theta Z}]}{m_i!}\nonumber\\
&\qquad\qquad
+
n!\mathrm{e}^{\theta}\sum_{i=1}^{n}\frac{G_{P}^{(i)}(\mathbb{E}[\mathrm{e}^{\theta Z}])}{i!}
\sum_{m_1+m_2+\ldots+m_i=n}\frac{
\frac{\mathrm{d}^{m_1}}{\mathrm{d}\theta^{m_1}}
\mathbb{E}[\mathrm{e}^{\theta Z}]}{m_1!}\ldots\frac{\frac{\mathrm{d}^{m_i}}{\mathrm{d}\theta^{m_i}}
\mathbb{E}[\mathrm{e}^{\theta Z}]}{m_i!}\nonumber\\
&=\mathrm {e}^\theta G_P(\mathbb E\mathrm{e}^{\theta Z})+
\mathrm{e}^{\theta}\sum_{k=1}^{n-1}\binom{n}{k}
k!\sum_{i=1}^{k}\frac{G_{P}^{(i)}(\mathbb{E}[\mathrm{e}^{\theta Z}])}{i!}
\sum_{m_1+m_2+\ldots+m_i=k}\frac{
\frac{\mathrm{d}^{m_1}}{\mathrm{d}\theta^{m_1}}
\mathbb{E}[\mathrm{e}^{\theta Z}]}{m_1!}\ldots\frac{\frac{\mathrm{d}^{m_i}}{\mathrm{d}\theta^{m_i}}
\mathbb{E}[\mathrm{e}^{\theta Z}]}{m_i!}\label{eq:02112023II}\\
&\qquad
+
n!\mathrm {e}^\theta\sum_{i=2}^{n}\frac{G_{P}^{(i)}(\mathbb{E}[\mathrm{e}^{\theta Z}])}{i!}
\sum_{m_1+m_2+\ldots+m_i=n}\frac{
\frac{\mathrm{d}^{m_1}}{\mathrm{d}\theta^{m_1}}
\mathbb{E}[\mathrm{e}^{\theta Z}]}{m_1!}\ldots\frac{\frac{\mathrm{d}^{m_i}}{\mathrm{d}\theta^{m_i}}
\mathbb{E}[\mathrm{e}^{\theta Z}]}{m_i!}+{\mathrm {e}^\theta}G_{P}^{'}(\mathbb{E}[\mathrm{e}^{\theta Z}])\frac{\mathrm{d}^n}{\mathrm{d}\theta^n}\mathbb{E}[\mathrm{e}^{\theta Z}].\nonumber
\end{align}
Therefore,
\begin{align}
&(1-\mathrm {e}^\theta G_{P}^{'}(\mathbb{E}[\mathrm{e}^{\theta Z}]))\frac{\mathrm{d}^n}{\mathrm{d}\theta^n}\mathbb{E}[\mathrm{e}^{\theta Z}]\nonumber\\
&\quad
=\mathrm {e}^\theta G_P(\mathbb E\mathrm{e}^{\theta Z})+
\mathrm{e}^{\theta}\sum_{k=1}^{n-1}\binom{n}{k}
k!\sum_{i=1}^{k}\frac{G_{P}^{(i)}(\mathbb{E}[\mathrm{e}^{\theta Z}])}{i!}
\sum_{m_1+m_2+\ldots+m_i=k}\frac{
\frac{\mathrm{d}^{m_1}}{\mathrm{d}\theta^{m_1}}
\mathbb{E}[\mathrm{e}^{\theta Z}]}{m_1!}\ldots\frac{\frac{\mathrm{d}^{m_i}}{\mathrm{d}\theta^{m_i}}
\mathbb{E}[\mathrm{e}^{\theta Z}]}{m_i!}\nonumber\\
&\qquad
+
n!\mathrm {e}^\theta \sum_{i=2}^{n}\frac{G_{P}^{(i)}(\mathbb{E}[\mathrm{e}^{\theta Z}])}{i!}
\sum_{m_1+m_2+\ldots+m_i=n}\frac{
\frac{\mathrm{d}^{m_1}}{\mathrm{d}\theta^{m_1}}
\mathbb{E}[\mathrm{e}^{\theta Z}]}{m_1!}\ldots\frac{\frac{\mathrm{d}^{m_i}}{\mathrm{d}\theta^{m_i}}
\mathbb{E}[\mathrm{e}^{\theta Z}]}{m_i!}.\nonumber
\end{align}
Letting $\theta\uparrow 0$ in this relation we have
\begin{align}
(1-G_{P}^{'}(1))\mathbb{E}[Z^n]&
=1+\sum_{k=1}^{n-1}\binom{n}{k}
k!\sum_{i=1}^{k}\frac{G_{P}^{(i)}(1)}{i!}
\sum_{m_1+m_2+\ldots+m_i=k}\frac{
\mathbb{E}[Z^{m_1}]}{m_1!}\ldots\frac{
\mathbb{E}[Z^{m_i}]}{m_i!}\nonumber\\
&\qquad
+
n!\sum_{i=2}^{n}\frac{G_{P}^{(i)}(1)}{i!}
\sum_{m_1+m_2+\ldots+m_i=n}\frac{\mathbb{E}[Z^{m_1}]}{m_1!}\ldots\frac{\mathbb{E}[Z^{m_i}]}{m_i!}.\nonumber
\end{align}
Since
\begin{equation}\label{eq:02112023III}
G_{P}^{'}(1)=\mathbb E(P)_1=\mathbb E P\quad\text{and}\quad G_{P}^{(i)}(1)=\mathbb{E}[P(P-1)\ldots(P-(i-1))],\quad 2\leq i\leq n,
\end{equation}
we have
\begin{align}
\mathbb{E}[Z^n]&
=\frac{1}{1- \mathbb E P}\Biggl(1+\sum_{k=1}^{n-1}\binom{n}{k}
k!\sum_{i=1}^{k}\frac{\mathbb E(P)_i}{i!}
\sum_{m_1+m_2+\ldots+m_i=k}\frac{
\mathbb{E}[Z^{m_1}]}{m_1!}\ldots\frac{
\mathbb{E}[Z^{m_i}]}{m_i!}\nonumber\\
&\qquad
+
n!\sum_{i=2}^{n}\frac{\mathbb E(P)_i}{i!}
\sum_{m_1+m_2+\ldots+m_i=n}\frac{\mathbb{E}[Z^{m_1}]}{m_1!}\ldots\frac{\mathbb{E}[Z^{m_i}]}{m_i!}\Biggr),\quad n\in\mathbb N.\label{eq:02112023I}
\end{align}

Reasoning by induction on $n\in\mathbb N$,  by relation \eqref{eq:02112023I}
we immediately have that $\mathbb E Z^n<\infty$ (note that, for $n=1$, we have 
$\mathbb E Z=1/(1-\mathbb E P)<\infty$,  and that, for $n\geq 2$,
all the moments of $Z$ involved in the right-hand side of \eqref{eq:02112023I} are of order less than or equal to $n-1$).

The formula \eqref{eq:recursive} follows by letting $\theta\uparrow 0$ in \eqref{eq:02112023II} and using the equalities in \eqref{eq:02112023III}.
\\
\noindent{\it Step\,\,3: Computing\,\,$\mathbb E Z^2$,\,\,$\mathbb E Z^3$\,\,and\,\,$\mathbb E Z^4$.} The claimed expressions for the second,  the third and the fourth moments of $Z$ easily follow by \eqref{eq:recursive} (or \eqref{eq:02112023I}).

For instance,  as far as the second moment is concerned,  the formula gives
\[
\mathbb E Z^2=1+2\mathbb E P\mathbb E Z+\mathbb E P\mathbb E Z^2+\mathbb E P(P-1)(\mathbb E Z)^2,
\]
from which the claimed expression of $\mathbb E Z^2$ immediately follows (recalling that $\mathbb E Z=(1-\mathbb E P)^{-1}$). We omit the computations for the third and the fourth moments of $Z$.
\\
\noindent$\square$

\subsection{Generalized compound Hawkes processes with Poisson offspring distribution}

\subsubsection{Gaussian approximation}\label{subsec:gausspoisson}

In this paragraph we suppose that $Z$ is distributed as the total progeny of a Galton-Watson process with one ancestor and offspring distribution the Poisson law with mean $h\in (0,1)$,  and that $\{X_n\}_{n\geq 1}$ is a Poisson process on $\R^d$ with intensity function $\lambda(x)=\lambda\bold{1}_B(x)$,  $x\in\mathbb R^d$,  for some $\lambda>0$ and some Borel set $B\subseteq\mathbb R^d$.  We denote by $V_{\mathrm{Poisson}}$ the corresponding generalized compound Hawkes process and by $W_{\mathrm{Poisson}}$ the functional \eqref{eq:11102023terzoV}
with $V_{\mathrm{Poisson}}$ in place of $V$.

\begin{Corollary}\label{cor:11102023quarto}
Under the foregoing assumptions and notation,  if the Borel sets $B$ and $C$ are such that $0<\mathrm{Leb}(B\cap C)<+\infty$ and $\E M^2\in (0,\infty)$,  then the bounds \eqref{eq:11102023primo} and \eqref{eq:11102023secondo} hold with $W_{\mathrm{Poisson}}$
in place of $W$,
\begin{equation}
\label{eq:05062023uno}
\mathbb E Z^3=\frac{1+2h}{(1-h)^5}
\end{equation}
and
\begin{align}
    \mathbb E Z^4=\frac{1}{1-h} \Bigg [1+ \frac{4h}{1-h} +\frac{6h^2}{(1-h)^2}+\frac{4h^3+6h}{(1-h)^3} + \frac{h^4+12h^2}{(1-h)^4 }
    +\frac{6h^3}{(1-h)^5}+\frac{11h^2+4h}{(1-h)^6}\Bigg].
\label{eq:05062023due}
\end{align} 
\end{Corollary}
\noindent{\it Proof}. \\ 
Note that conditions \eqref{eq:sub}  and \eqref{eq:2023} are satisfied (the latter holds for any $n\in\mathbb N$).  Therefore,  the expressions of the third and the fourth moment of the total progeny are given by Proposition \ref{thm:branching}.  The claim follows by Corollary \ref{cor:2}. 
\\
\noindent$\square$

\subsubsection{Moderate deviations,  Bernstein-type concentration inequalities and Normal approximation bounds with Cram\'er correction term}\label{sec:22112023}

In this paragraph we suppose that,  for each $\ell\in\mathbb N$,  $Z=Z_1^{\ell}(\R^d,\R)$ is distributed as the total progeny of a Galton-Watson process with one ancestor and offspring distribution the Poisson law with mean $h\in (0,1)$,  and that $\{X_n^{(\ell)}\}_{n\geq 1}$ is a Poisson process on $\R^d$ with intensity function $\lambda_\ell (x)=\lambda_\ell \boldsymbol 1_{B_\ell}(x)$,  $x\in\mathbb R^d$,  for positive constants $\lambda_\ell>0$ and Borel sets $B_\ell\subseteq\mathbb{R}^d$,  $\ell \in \mathbb N$.  We denote by $V_{\mathrm{Poisson}}^{(\ell)}$ the corresponding generalized compound Hawkes process and by $W_{\mathrm{Poisson}}^{(\ell)}$ the functional \eqref{eq:11102023terzoVelle} with $V_{\mathrm{Poisson}}^{(\ell)}$ in place of $V_\ell$.

\begin{Corollary}\label{cor:11102023quinto}
Let the foregoing assumptions and notation prevail,  and let the Borel sets $B_\ell$ and $C_\ell$, $\ell \in \mathbb N$,  be such that $0<\mathrm{Leb}(B_\ell \cap C_\ell)<+\infty$,  $\ell\in\mathbb N$,  $\E M^2>0$ and
\begin{equation}\label{eq:newassump}
\frac{\E |M|^m}{(\E M^2)^{m/2}}\leq (m!)^{\gamma},\quad\text{for any $m\geq 3$ and some $\gamma\geq 0$.}
\end{equation}
Then:\\
\noindent$(i)$ If $h-1-\log h\geq 1$,  then the sequence $\{W_{\mathrm{Poisson}}^{(\ell)}(C_\ell)\}_{\ell\geq 1}$ satisfies $\bold{MDP}(\gamma,\{\Delta_\ell\}_{\ell\in\mathbb N})$,  $\bold{BCI}(\gamma,\{\Delta_\ell\}_{\ell\in\mathbb N})$,  $\bold{NACC}(\gamma,\{\Delta_\ell\}_{\ell\in\mathbb N})$ where
\[
\Delta_\ell:=h\sqrt{\lambda_\ell\text{Leb}(B_\ell \cap C_\ell)}.
\]
\noindent$(ii)$ If $h-1-\log h<1$,  then the sequence $\{W_{\mathrm{Poisson}}^{(\ell)}(C_\ell)\}_{\ell\geq 1}$ satisfies $\bold{MDP}(\gamma,\{\Delta_\ell\}_{\ell\in\mathbb N})$,  $\bold{BCI}(\gamma,\{\Delta_\ell\}_{\ell\in\mathbb N})$,  $\bold{NACC}(\gamma,\{\Delta_\ell\}_{\ell\in\mathbb N})$ where
\[
\Delta_\ell := h(h-1-\log h)^3\sqrt{\lambda_\ell \text{Leb}(B_\ell \cap C_\ell)}.
\]
\end{Corollary}

\begin{Example}\label{re:examples}
Note that the condition \eqref{eq:newassump} is trivially satisfied with $\gamma=0$ if $M$ is a constant different from zero. Similarly, if $M$ were a uniform variable on $[0,D]$ with $D> 0$,  the moments of $M$ verify 
$$\frac{\E |M|^m}{(\E M^2)^{m/2}}=\frac{3^{\frac{m}{2}}}{m+1}<m!.$$
Assumption \eqref{eq:newassump} holds with $\gamma=1$ even if $M$ is exponentially distributed with mean $\mu^{-1}$, for some $\mu>0$.  Indeed, in such a case we have
\[
\frac{\E |M|^m}{(\E M^2)^{m/2}}=\frac{m!}{2^{m/2}}<m!.
\]
Another example is $M$ Gaussian distributed with mean zero and variance $\sigma^2$. Indeed,  in such a case we have 
$$\frac{\E |M|^m}{(\E M^2)^{m/2}} \leq (m-1)!!.$$
By distinguishing the two cases $m=2p+1$ and $m=2p+2$ for $p=2,3,\cdots$, we conclude that 
$$\frac{\E |M|^m}{(\E M^2)^{m/2}}\leq \sqrt {m!},$$
which entails that condition \eqref{eq:newassump} is satisfied with $\gamma=\frac{1}{2}.$
\end{Example}

The proof of the corollary exploits the following lemmas,  which are shown in Section \ref{sec:lemmas}.
\begin{Lemma}
    \label{lemma:abel-plana}
    For any $\nu>0$ and any integer $m\geq 2$ we have
    $$\sum_{k=1}^{+\infty}\mathrm{e}^{-\nu k}k^{m-1}=\nu ^{-m} (m-1)!+ R_m,$$
    with 
    $$|R_m| \leq \frac{1}{\pi(m-1)}+ \frac{2(m-1)!}{\pi^m} .$$
\end{Lemma}

\begin{Lemma}
\label{lemma:inequality}
The function
\[
f(x):=x(x-1-\log x)^2,\quad x\in (0,1)
\]
is such that $f(x)\in (0,1)$.

\end{Lemma}

\noindent{\it Proof\,\,of\,\,Corollary\,\,\ref{cor:11102023quinto}.}\\
It is well-known that the total progeny $Z$ of a sub-critical Galton-Watson process with one ancestor and Poisson offspring law with mean $h\in (0,1)$ follows the Borel distribution (\textit{cf.} \cite{P2} and the references therein), i.e.,  
$$\mathbb P (Z=k)=\frac{\mathrm e^{-hk}(hk)^{k-1}}{k!}, \quad k=1,2,\ldots$$
Therefore,  by Stirling's inequality,  for $m\geq 3$ we have
\begin{align*}
    \E Z^m = \sum_{k=1}^{+\infty} \frac{\mathrm e^{-hk}(hk)^{k-1}}{k!} k^m
    &\leq \sum_{k=1}^{+\infty}\mathrm e^{-hk}(hk)^{k-1}  k^m \left(\frac{\mathrm e}{k} \right)^k \frac{1}{\sqrt {2\pi k}}\\
    &= \frac{1}{\sqrt {2\pi}} \sum_{k=1}^{+\infty} \mathrm e^{(1-h)k} \frac{k^m}{k\sqrt k} h^{k-1}\\
    &= \frac{1}{h\sqrt {2\pi }} \sum_{k=1}^{+\infty} \mathrm e^{-(h-1-\log h)k} k^{m-1.5} \\
    &\leq \frac{1}{h\sqrt {2\pi }} \sum_{k=1}^{+\infty} \mathrm e^{-(h-1-\log h)k} k^{m-1}.
\end{align*}

Using Lemma \ref{lemma:abel-plana} with $\nu:=h-1-\log h>0$,  we have
\begin{equation}\label{eq:03112023I}
\E Z^m \leq \frac{(m-1)!}{h\sqrt{2\pi}}\left(\nu^{-m}+\frac{2}{\pi^m}+  \frac{1}{\pi (m-1)(m-1)!}\right),\quad\text{$m\geq 3$.}
\end{equation}
We continue distinguishing the proofs of Part $(i)$ and of Part $(ii)$, in both cases we shall apply Corollary \ref{cor:ST}.\\
\noindent{\it Proof\,\,of\,\,Part\,\,(i)}.
If $\nu\geq 1$,  then, for any $m\geq 3$,
\begin{equation}\label{eq:13102023primo}
\nu^{-m}+\frac{2}{\pi^m}+  \frac{1}{\pi (m-1)(m-1)!} \leq 3\leq m.
\end{equation}
Combining this inequality with \eqref{eq:03112023I},  for $m\geq 3$ and $\ell\in\mathbb N$,  we have
\begin{align*}
\frac{\mathbb E Z^m}{\sqrt{\lambda_\ell \text{Leb}(B_\ell \cap C_\ell)}^{m-2} }\leq& \frac{1}{h\sqrt{2\pi}}\frac{m!}{\sqrt{\lambda_\ell \text{Leb}(B_\ell \cap C_\ell)}^{m-2}}\\
=& \frac{h^{m-3}}{\sqrt{2\pi}}\frac{m!}{(h\sqrt{\lambda_\ell \text{Leb}(B_\ell \cap C_\ell)})^{m-2}}\\
\leq & \frac{m!}{(h\sqrt{\lambda_\ell \text{Leb}(B_\ell \cap C_\ell)})^{m-2}}
\end{align*}
since $h<1$ and $\sqrt{2\pi}>1$.  Combining this latter inequality with the assumption \eqref{eq:newassump} 
we have that condition \eqref{hyp:29092023} is satisfied with $\Delta_\ell=h\sqrt{\lambda_\ell \text{Leb}(B_\ell \cap C_\ell)}$,  and the claim follows by Corollary \ref{cor:ST}.\\
\noindent{\it Proof\,\,of\,\,Part\,\,(ii)}.  Now we suppose $\nu < 1$.  Since

\begin{equation}\label{eq:13102023secondo}
\frac{2}{\pi^m}+\frac{1}{\pi (m-1)(m-1)!}\leq\frac{2}{\pi}\left(\frac{1}{\pi^2}+\frac{1}{8}\right)<0.16,\quad\text{for any $m\geq 3$}
\end{equation}

by \eqref{eq:03112023I} we have
\begin{align*}
\frac{\mathbb E Z^m}{\sqrt{\lambda_\ell \text{Leb}(B_\ell \cap C_\ell)}^{m-2}} \leq & \frac{1}{h\sqrt{2\pi}}\frac{(m-1)!}{\sqrt{\lambda_\ell \text{Leb}(B_\ell \cap C_\ell)}^{m-2}}\left(\frac{1}{\nu^{m}}+0.16\right)\\
    = &\frac{(m-1)!}{(\sqrt{\lambda_\ell \text{Leb}(B_\ell \cap C_\ell)}\nu)^{m-2}}\frac{\nu^m}{h\nu^2\sqrt{2\pi}}\left(\frac{1}{\nu^{m}}+0.16\right),\quad\text{$m\geq 3$,  $\ell\in\mathbb N$}.\\ 
\end{align*}
Since $\nu <1$,  we have 
\begin{align*}
    \frac{\nu^m}{h \nu^2\sqrt{2\pi} }\left(\frac{1}{\nu^{m}}+0.16\right) &\leq\frac{1.16}{h\nu^2\sqrt{2\pi}}<\frac{1}{h\nu^2}:=u_h.
\end{align*}

Therefore
\begin{align*}
\frac{\mathbb E Z^m}{\sqrt{\lambda_\ell \text{Leb}(B_\ell \cap C_\ell )}^{m-2}} 
    \leq &\frac{(m-1)!}{(\sqrt{\lambda_\ell \text{Leb}(B_\ell \cap C_\ell)}\nu)^{m-2}}u_h\\ 
    =&\frac{(m-1)!}{(\sqrt{\lambda_\ell \text{Leb}(B_\ell \cap C_\ell)}\nu u_h^{-1})^{m-2}}u_h^{3-m},\quad m\geq 3, \ell\in\mathbb N.
\end{align*}
By Lemma \ref{lemma:inequality} we have $u_h^{-1} <1$. Therefore,  for any $m\geq 3$ and $\ell\in\mathbb N$,  we have
\begin{align*}
\frac{\mathbb E Z^m}{\sqrt{\lambda_\ell \text{Leb}(B_\ell \cap C_\ell)}^{m-2}} 
    \leq & \frac{m!}{(\sqrt{\lambda_\ell \text{Leb}(B_\ell \cap C_\ell)}\nu u_h^{-1})^{m-2}}.
\end{align*}
Combining this latter inequality with the assumption \eqref{eq:newassump} we have that
condition \eqref{hyp:29092023} is satisfied with 
$\Delta_\ell := h\nu^3\sqrt{\lambda_\ell \text{Leb}(B_\ell \cap C_\ell)}$, and the claim follows by Corollary \ref{cor:ST}. 
\\
\noindent$\square$

\subsection{On the Gaussian approximation bound in the Kolmogorov distance and the Normal approximation with Cram\'er correction term}

The aim of this section is to illustrate,  by means of a simple example,  the differences
and the analogies between Gaussian approximation bounds in the Kolmogorov distance and Normal approximations with Cram\'er correction term.

Let $W_{\mathrm{Poisson}}^{(\ell)}$,  $\ell\in\mathbb N$,  be defined as at the beginning of Section  \ref{sec:22112023} with $\{(B_\ell,C_\ell)\}_{\ell\in\mathbb N}\subset\mathcal{B}(\mathbb{R}^d)^2$ a sequence of Borel sets such that $0<\mathrm{Leb}(B_\ell\cap C_\ell)<\infty$,  $\ell\in\mathbb N$,  
and $\lambda_\ell\mathrm{Leb}(B_\ell\cap C_\ell)\to+\infty$, as $\ell\to+\infty$.  Assume $M\equiv 1$ (so that \eqref{eq:newassump} holds with $\gamma=0$),  and let $\Delta_\ell$ be defined as in Corollary \ref{cor:11102023quinto}$(i)$ or $(ii)$.  We know that
$\{W_{\mathrm{Poisson}}^{(\ell)}(C_\ell)\}_{\ell\geq 1}$ satisfies $\bold{NACC}(0,\{\Delta_\ell\}_{\ell\in\mathbb N})$.  

For a fixed $0<r<1/3$,  let $\ell^*$ be sufficiently large that
\[
0<\Delta_{\ell}^r<c_0\Delta_{\ell},\quad\text{for all $\ell\geq\ell^*$}
\]
where $c_0$ is the positive constant which appears in the definition of the $\bold{NACC}(0,\{\Delta_\ell\}_{\ell\in\mathbb N})$.  Setting $x_\ell:=\Delta_{\ell}^r$,  for all $\ell\geq\ell^*$,  we have
\begin{align}
\mathbb{P}(W_{\mathrm{Poisson}}^{(\ell)}(C_\ell)\geq x_\ell)-\mathbb{P}(G\geq x_\ell)&=\mathrm{e}^{L_{\ell,x_\ell}^+}\mathbb P(G\geq x_\ell)\left(1+c_1\theta_{\ell,x_\ell}^+\frac{1+x_\ell}{\Delta_{\ell}}\right)-\mathbb{P}(G\geq x_\ell)\nonumber\\
&=(\mathrm{e}^{L_{\ell,x_\ell}^+}-1)\mathbb{P}(G\geq x_\ell)+
c_1\theta_{\ell,x_\ell}^{+}\mathrm{e}^{L_{\ell,x_\ell}^+}\mathbb P(G\geq x_\ell)\frac{1+x_\ell}{\Delta_{\ell}}.\nonumber
\end{align}
Since $|L_{\ell,x_\ell}^{+}|\leq c_2\Delta_\ell^{3r-1}$ and $|\theta_{\ell,x_\ell}^+|\leq 1$,  we have
\begin{align}
|\mathbb{P}(W_{\mathrm{Poisson}}^{(\ell)}(C_\ell)\geq x_\ell)-\mathbb{P}(G\geq x_\ell)|
&\leq \left(c_2\Big|\frac{\mathrm{e}^{L_{\ell,x_\ell}^+}-1}{L_{\ell,x_\ell}^+}\Big|\Delta_\ell^{3r-1}+
c_1 O(1)\Big|\frac{1+x_\ell}{\Delta_{\ell}}\Big|\right) \mathbb P(G\geq x_\ell)\nonumber\\
&=\left(c_2 O(1)\Delta_\ell^{3r-1}+
c_1 O(1)\Big|\frac{1+\Delta^r_\ell}{\Delta_{\ell}}\Big|\right) \mathbb P(G\geq x_\ell)\nonumber\\
&=O(\Delta_\ell^{3r-1})\mathbb P(G\geq x_\ell).
\nonumber
\end{align}
Bounding the Gaussian tail from above we obtain
\begin{align}
    |\mathbb{P}(W_{\mathrm{Poisson}}^{(\ell)}(C_\ell)\geq x_\ell)-\mathbb{P}(G\geq x_\ell)|&\leq O(\Delta_\ell^{3r-1})\frac{\mathrm{e}^{-x_\ell^2/2}}{x_\ell \sqrt{2\pi}}\nonumber\\
    &=O(\Delta_\ell^{2r-1}\mathrm{e}^{-\Delta_\ell ^{2r}/2})\nonumber\\
    &=O\left(\frac{\exp\left({-\frac{u_h^{2r}(\lambda_\ell\mathrm{Leb}(B_\ell\cap C_\ell))^r}{2}}\right)}{(\lambda_\ell\mathrm{Leb}(B_\ell\cap C_\ell))^{\frac{1}{2}-r}} \right),\quad\text{as $\ell\to+\infty$}\label{eq:14112023primo}
\end{align}
where either $u_h:=h$ or $u_h:=h(h-1-\log h)^3$.  On the other hand,  the bound \eqref{eq:11102023secondo} and the relations \eqref{eq:05062023uno} and \eqref{eq:05062023due} yield
\begin{align}
|\mathbb{P}(W_{\mathrm{Poisson}}^{(\ell)}(C_\ell)\geq x_\ell)-\mathbb{P}(G\geq x_\ell)|
=O\left(\frac{1}{(\lambda_\ell\mathrm{Leb}(B_\ell\cap C_\ell))^{1/2}}\right),\quad\text{as $\ell \to +\infty$ .}\label{eq:14112023secondo}
\end{align}
Clearly the rate \eqref{eq:14112023primo} is much faster than \eqref{eq:14112023secondo}.  

Now let $x\in (0,\infty)$ be arbitrarily fixed.  For $\ell$ large enough we have
$x\in [0,c_0\Delta_\ell]$,  and so by Corollary \ref{cor:11102023quinto}, for all $\ell$ large enough,  we have
\begin{align}
|\mathbb{P}(W_{\mathrm{Poisson}}^{(\ell)}(C_\ell)\geq x)-\mathbb{P}(G\geq x)|&\leq |\mathrm{e}^{L_{\ell,x}^+}-1|\mathbb{P}(G\geq x)+c_1\mathrm{e}^{L_{\ell,x}^+}\mathbb P(G\geq x)\theta_{\ell,x}^+\frac{1+x}{\Delta_\ell}\nonumber\\
&\leq c_2\Big|\frac{\mathrm{e}^{L_{\ell,x}^+}-1}{L_{\ell,x}^+}\Big|\frac{x^3}{\Delta_\ell}\mathbb{P}(G\geq x)+c_1\mathrm{e}^{L_{\ell,x}^+}\mathbb P(G\geq x)\theta_{\ell,x}^+\frac{1+x}{\Delta_\ell}\nonumber\\
&=O\left(\frac{1}{\sqrt{\lambda_\ell\mathrm{Leb}(B_\ell\cap C_\ell)}}\right).\nonumber
\end{align}
Clearly,  the same rate is provided by the bound \eqref{eq:11102023secondo} and the relations \eqref{eq:05062023uno} and \eqref{eq:05062023due}.  We emphasize that: $(i)$ The inequality \eqref{eq:11102023secondo} and the relations \eqref{eq:05062023uno} and \eqref{eq:05062023due}
yield indeed an explicit bound on the quantity $|\mathbb{P}(W_{\mathrm{Poisson}}^{(\ell)}(C_\ell)\geq x)-\mathbb{P}(G\geq x)|$,  for any $x\in\mathbb R$ and any $\ell\in\mathbb N$,  $(ii)$ An explicit bound on the quantity $|\mathbb{P}(W_{\mathrm{Poisson}}^{(\ell)}(C_\ell)\geq x)-\mathbb{P}(G\geq x)|$,  for any $x\in\mathbb R$ and any $\ell\in\mathbb N$,  is not amenable via the Normal approximation with Cram\'er correction term (for various obvious reasons).

\subsection{Comparison with some related literature}

\subsubsection{Gaussian approximation}\label{subsec:21122023primo}

Let $N$ be a classical Hawkes process on $(0,\infty)$ with parameters $(\lambda,g)$.
Then Corollary \ref{cor:11102023quarto} with $B=(0,\infty)$ and $C=(0,t]$,  $t>0$, gives {\em explicit bounds} for the Gaussian approximation of 
\[
\frac{N((0,t])-\E N((0,t])}{\sqrt{\mathbb V\mathrm{ar}(N((0,t]))}}, 
\]
both in the Wasserstein and Kolmogorov distances.  Note that on the fertility function $g$ we only assume the standard stability condition $h:=\int_0^\infty g(t)\,\mathrm{d}t\in (0,1)$.

It is worthwhile to compare these bounds with the ones in \cite{HHKR} and \cite{KPR}.  

Theorem 3.13 in \cite{HHKR} gives a bound of the kind
\[
d_W\left(\frac{N((0,t])-\mathbb E N((0,t])}{\sqrt{t}},G\right)\leq c/\sqrt t, \quad t>0
\]
for some constant $c>0$ which is not explicitly computed and for specific choices of $g$ (exponential and Erlang). This result has been extended in \cite{KPR} to fertility functions $g: [0,\infty)\to [0,\infty)$ such that $h\in (0,1)$ and $\int_0^\infty t g(t)\mathrm{d}t<\infty$. The techniques used in \cite{HHKR} and \cite{KPR} are based on the Poisson embedding construction of Hawkes processes and the Malliavin calculus on the Poisson space.  These ideas were previously used in \cite{T1} and \cite{T2} for the purpose of Gaussian and Poisson approximation of some classes of nonlinear Hawkes processes.  Note that, in contrast with classical Hawkes processes,  nonlinear Hawkes processes (introduced in \cite{BM}) do not have a Poisson cluster representation.  For Poisson cluster processes (such as classical Hawkes processes),  the number of points on some measurable set can be represented as an integral with respect to a suitable Poisson random measure. 
  
As a consequence,  results on the Gaussian approximation of the number of points on a measurable set can be obtained by applying the general results in \cite{LPS}.   

\subsubsection{Moderate deviations}\label{subsec:21122023secondo}

In this section we compare Corollary \ref{cor:11102023quinto} with a couple of related
results in literature.

Firsty,  we prove that Corollary \ref{cor:11102023quinto},  when specialized to a classical Hawkes process on $(0,\infty)$,  implies the same moderate deviation principle provided in \cite{ZHU} (see Theorem 1 therein),  with an alternate assumption on the fertility function of the process; we refer the reader to \cite{GW} for sample-path moderate deviation principles,  on the space of càdlàg functions on $[0,1]$ equipped with the Skorokhod topology,  of Poisson cluster point processes on the line, see Theorem 2 therein.

Secondly,  we compare Corollary \ref{cor:11102023quinto} (again specialized to a classical Hawkes process on $(0,\infty)$) with Theorem 8 in \cite{GZ}.

Hereafter,  $N$ denotes a classical Hawkes process on $(0,\infty)$ with parameters $(\lambda,g)$ and we
assume that $g$ satisfies the usual stability condition $h:=\int_0^\infty g(t)\,\mathrm{d}t\in (0,1)$. 

Corollary \ref{cor:11102023quinto} with $\lambda_\ell=\lambda>0$,  $B_\ell=(0,\infty)$ and $C_\ell=(0,\ell]$,  $\ell\in\mathbb N$,  and $M\equiv 1$ yields that,
for any sequence of positive numbers $\{a_\ell\}_{\ell\in\mathbb N}$ such that $\lim_{\ell\to\infty}a_\ell=+\infty$ and $\lim_{\ell\to\infty}\frac{a_{\ell}}{\sqrt{\ell}}=0$,  for any Borel set $B\subseteq\R$,
\begin{align}
-\inf_{x\in\overset{\circ}B}\frac{x^2}{2}&\leq\liminf_{\ell\to\infty}a_\ell^{-2}\log\mathbb{P}\Biggl(
\frac{N((0,\ell])-\mathbb E N((0,\ell])}{a_\ell\sqrt{\mathbb{V}\mathrm{ar}(N((0,\ell]))}}
\in B\Biggr)\nonumber\\
&\leq
\limsup_{\ell\to\infty}a_\ell^{-2}\log\mathbb{P}\Biggl(
\frac{N((0,\ell])-\mathbb E N((0,\ell])}{a_\ell\sqrt{\mathbb{V}\mathrm{ar}(N((0,\ell]))}}
\in B\Biggr)
\leq-\inf_{x\in\overline B}\frac{x^2}{2}.\nonumber
\end{align}
Reasoning by contradiction,  one has that,  for any function $a(\cdot)$ such that
$\lim_{t\to\infty}a(t)=+\infty$ and $a(t)=o(\sqrt t)$, as $t\to+\infty$,  and any Borel set $B\subseteq\mathbb R$,  we have
\begin{align}
-\inf_{x\in\overset{\circ}B}\frac{x^2}{2}&\leq\liminf_{t\to\infty}a(t)^{-2}\log\mathbb{P}\Biggl(
\frac{N((0,t])-\mathbb E N((0,t])}{a(t)\sqrt{\mathbb{V}\mathrm{ar}(N((0,t]))}}
\in B\Biggr)\nonumber\\
&\leq
\limsup_{\ell\to\infty}a(t)^{-2}\log\mathbb{P}\Biggl(
\frac{N((0,t])-\mathbb E N((0,t])}{a(t)\sqrt{\mathbb{V}\mathrm{ar}(N((0,t]))}}
\in B\Biggr)
\leq-\inf_{x\in\overline B}\frac{x^2}{2}.\nonumber
\end{align}

It is easily realized that (\textit{cf.} \cite{Bacry} for example),  under the stability condition,
\[
\mathbb E N((0,t])/t\to\frac{\lambda}{1-h}\quad\text{and}\quad\mathbb{V}\mathrm{ar}(N((0,t]))/t\to\frac{\lambda}{(1-h)^3},\quad\text{as $t\to+\infty$.}
\]

So, letting $b(\cdot)$ denote a function such that $\sqrt t=o(b(t))$ and $b(t)=o(t)$,  as $t\to\infty$,
and setting $a(t):=b(t)/\sqrt{\mathbb{V}\mathrm{ar}(N((0,t]))}$, we have
that,  for any Borel set $B\subseteq\mathbb R$,  
\begin{align}
-\inf_{x\in\overset{\circ}B}\frac{x^2(1-h)^3}{2\lambda}&\leq\liminf_{t\to\infty}\frac{t}{b(t)^2}\log\mathbb{P}\Biggl(
\frac{N((0,t])-\mathbb E N((0,t])}{b(t)}
\in B\Biggr)\nonumber\\
&\leq
\limsup_{\ell\to\infty}\frac{t}{b(t)^2}\log\mathbb{P}\Biggl(
\frac{N((0,t])-\mathbb E N((0,t])}{b(t)}\in B\Biggr)
\leq-\inf_{x\in\overline B}\frac{x^2(1-h)^3}{2\lambda}.\label{eq:27112023secondo}
\end{align}
By Lemma 5 in \cite{Bacry},  we have that,  if in addition to the stability condition $h\in (0,1)$ we assume
\begin{equation}\label{eq:27112023prima}
\int_0^\infty\sqrt t g(t)\,\mathrm{d}t<\infty,
\end{equation}
then
\[
\frac{\mathbb E N((0,t])-\lambda t/(1-h)}{b(t)}\to 0,\quad\text{as $t\to\infty$.}
\]
So, for an arbitrarily fixed  $\delta>0$, there exists $t_\delta$ such that for any $t>t_\delta$ it holds
\[
\Big|\frac{N((0,t])-\mathbb E N((0,t])}{b(t)}-\frac{N((0,t])-\frac{\lambda t}{1-h}}{b(t)}\Big|=\Big|\frac{\mathbb E N((0,t])-\lambda t/(1-h)}{b(t)}\Big|<\delta.
\]
Therefore,  for an arbitrarily fixed  $\delta>0$, there exists $t_\delta$ such that for any $t>t_\delta$ we have
\[
\log\mathbb{P}\Biggl(\Big|\frac{N((0,t])-\mathbb E N((0,t])}{b(t)}-\frac{N((0,t])-\frac{\lambda t}{1-h}}{b(t)}\Big|
>\delta\Biggr)=-\infty.
\]
Hence the processes $\{\frac{N((0,t])-\mathbb E N((0,t])}{b(t)}\}_{t>0}$ and
$\{\frac{N((0,t])-\lambda t/(1-h)}{b(t)}\}_{t>0}$ are exponentially equivalent (see Definition 4.2.10 p. 130 in \cite{DZ})  and so by Theorem 4.2.13 p. 130 in \cite{DZ} relation \eqref{eq:27112023secondo} holds with $\mathbb E N((0,t])$ replaced by $\lambda t/(1-h)$.

Thus we recover the moderate deviation principle proved in \cite{ZHU} under an alternate condition on $g$ (in Zhu's paper it is assumed the stability condition and $\sup_{t>0}t^{3/2}g(t)<\infty$,  which is clearly different from \eqref{eq:27112023prima}).

Theorem 8 in \cite{GZ} states that,  under the assumption 
$$\int_0 ^{\infty} t g(t) \d t <+\infty,$$
(which is clearly stronger than  \eqref{eq:27112023prima}) for any $y(t)=o(t^{1/2-1/m})$,  as $t\to+\infty$,  $m\geq 3$ integer,   and positive function $b(\cdot)$,  it holds
\begin{equation*}
%    \label{eq:GZ}
    \mathbb P \left (\frac{N((0,t])-\frac{\lambda t}{1-h}}{b(t)} \geq \frac{\sqrt{t}y(t)}{b(t)} \frac{\sqrt{\lambda}}{(1-h)^{3/2}} \right)=\frac{(1+o(1))}{y(t) \sqrt{2\pi}} \mathrm e ^{-\sum_{i=2}^{m-1}c_i \frac{y(t)^i}{t^{(i-2)/2}}},\quad\text{as $t\to+\infty$}
\end{equation*}
where $\{c_i\}_{i=1,\cdots,m-2}$ are real coefficients that can be computed explicitly, for instance,  one has $c_2=\frac 12 $.  

In particular,  if 
$b(t)=o( t^{2/3})$,  as $t\to+\infty$,  then by choosing 
\[
y(t)=K \frac{(1-h)^{3/2}}{\sqrt {\lambda}}\frac{b(t)}{\sqrt{t}} =o(t^{1/2-1/3}),\quad\text{as $t\to+\infty$,  for some $K>0$}
\]
we have
\begin{equation}
\label{eq:GZ2}
    \mathbb P \left( \frac{N((0,t])-\frac{\lambda t}{1-h}}{b(t)} \geq K\right)= \frac{\sqrt \lambda (1+o(1)) }{(1-h)^{3/2}K}\frac{\sqrt t}{b(t)} \mathrm e ^{-\frac{K^2}{2} \frac{(1-h)^{3}}{\lambda}\frac{b(t)^2}{t}},
\end{equation}
which is a more precise form of the relation \eqref{eq:27112023secondo} for the Borel set $B=[K,+\infty)$. Note that,  unlike formula \eqref{eq:27112023secondo},  which is valid for any Borel set $B$,  formula \eqref{eq:GZ2} gives asymptotic estimates only for half lines. 
%\\
%\textcolor{red}{What about to cut this two lines?: The precise moderate deviation formula also holds for $t^{2/3} =o(b(t))$, albeit with extra terms as shown in Remark 11 in \cite{GZ}. Note however that their result is logarithmically equivalent to ours as $t\to \infty$. }

\subsubsection{Bernstein-type concentration inequalities}

In this section we present some consequences of Corollary \ref{cor:11102023quinto} concerning stationary compound Hawkes processes on the line,  i.e., $B_T:=\R$,  observed on the time interval $C_T:=(0,T]$; here $T$ replaces $\ell$ to emphasize the dependence on time. 

If we interpret the mark $M$ as the claim that an insurer must pay to an insurance policy holder,  the variable $V((0,T])$ (defined by \eqref{eq:rapprS}) represents the total loss incurred by the insurer in the time interval $(0,T]$. 

Assume that the claims arrivals are modeled by the points of a Hawkes process of baseline intensity $\lambda>0$ and Poisson offspring distribution with mean $h\in (0,1)$ satisfying $h-1-\log h \geq 1$.  Assume moreover that the mark $M$ follows the exponential distribution of parameter $\mu^{-1} \in (0,+\infty)$. Then, Corollary \ref{cor:11102023quinto} yields
 \[
\mathbb P \left( \left | \frac{V((0,T])-\E V((0,T])}{\sqrt{\mathbb{V}\mathrm{ar} V((0,T])}}\right | \geq x  \right) \leq 2 \exp \left (-\frac{1}{4} \min \left \{\frac{x^2}{2^{1+\gamma}}, (x\Delta _T)^{\frac{1}{1+\gamma}} \right \}\right),\quad\text{$x\geq 0$}
\]
where $\Delta_T=h\sqrt{\lambda T}$ and $\gamma=1$ by virtue of Example \ref{re:examples}.  By stationarity,  this inequality can be rewritten as
%rearranging the terms in the last inequality and using the fact that the loss process has reached its stationary regime, the last inequality can be put under the form 
\begin{align}
&\mathbb P \left (V((0,T]) \in  \left[\frac{\lambda \mu}{1-h} T-x \sqrt{\frac{2\mu^2 \lambda}{(1-h)^3}T }, \frac{\lambda \mu}{1-h} T+x \sqrt{\frac{2\mu^2 \lambda}{(1-h)^3}T }  \right]\right)\nonumber\\
&\qquad\qquad
\geq 1-2 \exp \left (-\frac{1}{4} \min \left \{\frac{x^2}{4}, (xh\sqrt{\lambda T})^{1/2} \right \}\right),\quad\text{$x\geq 0$}\nonumber
\end{align}
which yields a non-asymptotic lower bound on the probability that the total loss is within $x$ times its standard deviation.

%For more results on the Hawkes process relating to insurance in the Markov framework, we refer the interested reader to \cite{JD} and \cite{SZZ}.\\
 Another quantity of interest for insurers is the probability that the total loss exceeds greatly its expected value.  The Bernstein-type concentration inequality,  being valid for any $x\geq 0$, yields an upper bound on this probability.  Indeed,  by choosing 
\[
x=x_T=(k-1)\sqrt{\frac{\lambda(1-h)T}{2}},\quad\text{for some $k>1$}
\]
we have
 \begin{align}
     \mathbb P \left(V((0,T]) \geq d \frac{\lambda \mu}{1-h} T \right)&= \mathbb P \left (\frac{V((0,T]) - \E V((0,T])}{\sqrt{\mathbb{V}\mathrm{ar} V((0,T])}} \geq x_T\right)
\nonumber\\
     &\leq \mathbb P \left ( \left |\frac{V((0,T]) - \E V((0,T])}{\sqrt{\mathbb{V}\mathrm{ar} V((0,T])}}\right | \geq x_T\right)
\nonumber\\
     &\leq 2 \exp \left (-\frac{1}{4} \min \left \{\frac{(k-1)^2\lambda (1-h)T}{8}, \left((k-1)\lambda h \sqrt{\frac{1-h}{2}} T\right)^{1/2}\right \} \right).\label{eq:29052024I}
 \end{align}
If the time horizon $T$ satisfies $T \geq \frac{2^{11/2}h}{(d-1)^3 \lambda (1-h)^{3/2}}$,  then the inequality \eqref{eq:29052024I} simplifies to 
 \begin{equation}
 \label{ineq: Applied_BCI}\mathbb P \left(V((0,T]) \geq k \frac{\lambda \mu}{1-h} T \right) \leq 2\exp \left (-\frac{1}{4} \left((k-1)\lambda h \sqrt{\frac{1-h}{2}} T\right)^{1/2} \right). 
 \end{equation}
A similar (non-asymptotic) inequality appears in Proposition 2.1 of \cite{RR},  albeit working only for stationary Hawkes processes on the line whose fertility functions have a compact support,  and involving quantities that are not explicitly known. We also point out that, specializing inequality \eqref{ineq: Applied_BCI} for the simple Hawkes process, that is with constant marks ($\gamma=0$ by virtue of Example \ref{re:examples}), we find a decay rate for the tail probability similar to the one given in \cite{RR}, but with explicit constants.

\subsection{Generalized compound Hawkes processes with Binomial offspring distribution}

\subsubsection{Gaussian approximation}\label{subsec:gaussbinom}

In this paragraph we suppose that $Z$ is distributed as the total progeny of a Galton-Watson process with one ancestor and offspring distribution the Binomial law with parameters $(h,p)$,  with $h\in\mathbb N$ and $p\in (0,1)$ such that $hp\in (0,1)$.  We assume that $\{X_n\}_{n\geq 1}$ is a Poisson process on $\R^d$ with intensity function $\lambda(x)=\lambda\bold{1}_B(x)$,  $x\in\mathbb R^d$, for some positive constant $\lambda>0$ and some Borel set $B\subseteq\mathbb{R}^d$.  We denote by $V_{\mathrm{Binomial}}$ the corresponding generalized compound Hawkes process and by $W_{\mathrm{Binomial}}$ the functional 
\eqref{eq:11102023terzoV}
with $V_{\mathrm{Binomial}}$ in place of $V$.

\begin{Corollary}\label{cor:11102023quartoBIS}
Under the foregoing assumptions and notation,  if the Borel sets $B$ and $C$ are such that $0<\mathrm{Leb}(B\cap C)<+\infty$ and $\E M^2\in (0,\infty)$,  then the bounds \eqref{eq:11102023primo} and \eqref{eq:11102023secondo} hold with $W_{\mathrm{Binomial}}$
in place of $W$,
\begin{align}
&\mathbb E Z^3=\frac{1}{1-hp}\Biggl(1+\frac{3hp}{1-hp}+3\frac{p^2(h)_2}{(1-hp)^2}+\frac{p^3(h)_3+3hp(1-p)}{(1-hp)^3}
+3\frac{h^2p^2(1-p)^2}{(1-hp)^4}\Biggr)\nonumber
\end{align}
and
\begin{align*}
    \mathbb E Z^4=&\frac{1}{1-hp} \Bigg [1+\frac{4hp}{1-hp}+\frac{6p^2(h)_2}{(1-hp)^2}+4\frac{p^3(h)_3}{(1-hp)^3} + \frac{p^4(h)_4}{(1-hp)^4 }\\
    &+\frac{3(1-hp^2)}{(1-hp)^3}\Big(2hp+\frac{4p^2(h)_2}{1-hp}+p^2(h)_2\frac{1-hp^2}{(1-hp)^3}+\frac{2p^3(h)_3}{(1-hp)^2}\Big)+4\mathbb E Z^3 \frac{hp(1-p)}{1-hp}\Bigg].
\end{align*}
Here $(h)_n:=h(h-1)\ldots (h-(n-1))\bold{1}_{\{h\geq n\}}$
\end{Corollary}
\noindent{\it Proof}. \\ 
Similar to the proof of Corollary \ref{cor:11102023quarto}.\\
\noindent$\square$

\subsubsection{Moderate deviations,  Bernstein-type concentration inequalities and Normal approximation bounds with Cram\'er correction term}

In this paragraph we suppose that,  for each $\ell\in\mathbb N$,  $Z=Z_1^{\ell}(\R^d,\R)$ is distributed as the total progeny of a Galton-Watson process with one ancestor and offspring distribution the Binomial law with parameters $(h,p)$,  with $h\in\mathbb N$ and $p\in (0,1)$ such that $hp\in (0,1)$.  We assume that $\{X_n^{(\ell)}\}_{n\geq 1}$ is a Poisson process on $\R^d$  with intensity function $\lambda_\ell (x)=\lambda_\ell \boldsymbol 1_{B_\ell}(x)$,  $x\in\mathbb R^d$,  for positive constants $\lambda_\ell>0$ and Borel sets $B_\ell\subseteq\mathbb R^d$,  $\ell \in \mathbb N$.  We denote by $V_{\mathrm{Binomial}}^{(\ell)}$ the corresponding generalized compound Hawkes process and by $W_{\mathrm{Binomial}}^{(\ell)}$ the functional \eqref{eq:11102023terzoVelle} with $V_{\mathrm{Binomial}}^{(\ell)}$ in place of $V_\ell$.

\begin{Corollary}\label{cor:11102023nono}
Let the foregoing assumptions and notation prevail,  and let the Borel sets $B_\ell$ and $C_\ell$, $\ell \in \mathbb N$,  be such that $0<\mathrm{Leb}(B_\ell \cap C_\ell)<+\infty$,  $\ell\in\mathbb N$, 
$\E M^2>0$ and assume \eqref{eq:newassump}.  Then:\\

\noindent$(i)1$ If $h=1$ and $p\leq\mathrm{e}^{-1}$,  then the sequence $\{W_{\mathrm{Binomial}}^{(\ell)}(C_\ell)\}_{\ell\geq 1}$ satisfies $\bold{MDP}(\gamma,\{\Delta_\ell\}_{\ell\in\mathbb N})$,  $\bold{BCI}(\gamma,\{\Delta_\ell\}_{\ell\in\mathbb N})$,  $\bold{NACC}(\gamma,\{\Delta_\ell\}_{\ell\in\mathbb N})$ where
\[
\Delta_\ell:= \frac{p}{1.05(1-p)}\sqrt{\lambda_\ell \text{Leb}(B_\ell \cap C_\ell)}.
\]
\noindent$(i)2$ If $h=1$ and $p>\mathrm{e}^{-1}$,  then the sequence $\{W_{\mathrm{Binomial}}^{(\ell)}(C_\ell)\}_{\ell\geq 1}$ satisfies $\bold{MDP}(\gamma,\{\Delta_\ell\}_{\ell\in\mathbb N})$,  $\bold{BCI}(\gamma,\{\Delta_\ell\}_{\ell\in\mathbb N})$,  $\bold{NACC}(\gamma,\{\Delta_\ell\}_{\ell\in\mathbb N})$ where
\[
\Delta_\ell:=\frac{p(\log p)^4}{1.05(1-p)}\sqrt{\lambda_\ell \text{Leb}(B_\ell \cap C_\ell)}.
\]
\noindent$(ii)1$ If $h\geq 2$ and $ph(h(1-p)/(h-1))^{h-1}\leq\mathrm{e}^{-1}$,
then the sequence $\{W_{\mathrm{Binomial}}^{(\ell)}(C_\ell)\}_{\ell\geq 1}$ satisfies $\bold{MDP}(\gamma,\{\Delta_\ell\}_{\ell\in\mathbb N})$,  $\bold{BCI}(\gamma,\{\Delta_\ell\}_{\ell\in\mathbb N})$,  $\bold{NACC}(\gamma,\{\Delta_\ell\}_{\ell\in\mathbb N})$ where
\[
\Delta_\ell:=\left(1+\sqrt{1+\frac{1}{h-1}} \frac{\mathrm e^{\frac{1}{24\cdot25}}(1-p)}{p(h-1)\sqrt {2\pi}}\right)^{-1}\sqrt{\lambda_\ell \text{Leb}(B_\ell \cap C_\ell)}
\]
\noindent$(ii)2$ If $h\geq 2$ and $ph(h(1-p)/(h-1))^{h-1}>\mathrm{e}^{-1}$,
then the sequence $\{W_{\mathrm{Binomial}}^{(\ell)}(C_\ell)\}_{\ell\geq 1}$ satisfies $\bold{MDP}(\gamma,\{\Delta_\ell\}_{\ell\in\mathbb N})$,  $\bold{BCI}(\gamma,\{\Delta_\ell\}_{\ell\in\mathbb N})$,  $\bold{NACC}(\gamma,\{\Delta_\ell\}_{\ell\in\mathbb N})$ where
\[
\Delta_\ell :=-\frac{\left(1+\sqrt{1+\frac{1}{h-1}} \frac{\mathrm e^{\frac{1}{24\cdot25}}(1-p)}{p(h-1)\sqrt {2\pi}}\right)^{-1}\left[\log\left(ph\left( \frac{h(1-p)}{h-1}\right) ^{h-1}\right)\right]^3}{1.16}\sqrt{\lambda_\ell \text{Leb}(B_\ell \cap C_\ell)}.
\]
\end{Corollary}

\noindent{\it Proof.}\\
It is well-known that the total progeny $Z$ of a sub-critical Galton-Watson process with one ancestor and Binomial offspring law with parameters $(h,p)$ follows the 
Consul distribution, i.e.,  
\begin{equation}\label{eq:03112023IV}
\mathbb P (Z=k)=\frac{1}{k}\binom{kh}{k-1}p^{k-1}(1-p)^{k(h-1)+1}, \quad k=1,2,\ldots.
\end{equation}
By Stirling's upper and lower bounds on the factorial,  for $k\geq 1$ and $h\geq 2$,  we have
\begin{align}
\binom{kh}{k-1}&=\frac{(kh)!}{(k-1)!(k(h-1)+1)!}\nonumber\\
&\leq\frac{1}{(k-1)!}\times\frac{\sqrt{2\pi kh}\left(\frac{kh}{\mathrm e}\right)^{kh}\mathrm{e}^{\frac{1}{12 kh}}}{\sqrt{2\pi (k(h-1)+1)}\left(\frac{k(h-1)+1}{\mathrm e}\right)^{k(h-1)+1}\mathrm{e}^{\frac{1}{12[k(h-1)+1]+1}}}\nonumber\\
&=\frac{\mathrm e^{-k+1}}{(k-1)!}\sqrt{\frac{kh}{k(h-1)+1}}\frac{(kh)^{kh}
}{(k(h-1)+1)^{k(h-1)+1}}\mathrm{e}^{\frac{1}{12 kh}-\frac{1}{12[k(h-1)+1]+1}}\nonumber\\
&=\frac{\mathrm e^{-k+1}}{(k-1)!}\sqrt{\frac{kh}{k(h-1)+1}}\frac{(k(h-1))^{kh}
}{(k(h-1)+1)^{k(h-1)+1}}\left(\frac{h}{h-1}\right)^{kh}\mathrm{e}^{\frac{1}{12 kh}-\frac{1}{12[k(h-1)+1]+1}}\nonumber\\
&\leq\frac{\mathrm e^{-k+1}}{(k-1)!}\sqrt{1+\frac{1}{h-1}}\frac{(k(h-1))^{kh}
}{(k(h-1)+1)^{k(h-1)+1}}\left(\frac{h}{h-1}\right)^{kh}\mathrm{e}^{\frac{1}{12 kh}-\frac{1}{12[k(h-1)+1]+1}}\label{eq:09102023}\\
&=\frac{\mathrm e^{-k+1}}{(k-1)!}\sqrt{1+\frac{1}{h-1}}\left(\frac{k(h-1)}{k(h-1)+1}\right)^{k(h-1)+1}[k(h-1)]^{k-1}
\left(\frac{h}{h-1}\right)^{kh}\mathrm{e}^{\frac{1}{12 kh}-\frac{1}{12[k(h-1)+1]+1}}\nonumber\\
&\leq\frac{\mathrm{e}^{\frac{1}{24\cdot 25}}}{h-1}\frac{\mathrm e^{-k+1} k^{k-1}}{(k-1)!}
\sqrt{1+\frac{1}{h-1}}\mathrm e ^{-1}
\left(\frac{h^h}{(h-1)^{h-1}}\right)^{k},\nonumber
\end{align}
where the inequality \eqref{eq:09102023} follows noticing that $kh/[k(h-1)+1]\leq h/(h-1)$ and the latter inequality follows noticing that
\[
\mathrm{e}^{\frac{1}{12 kh}-\frac{1}{12[k(h-1)+1]+1}}\leq\mathrm{e}^{\frac{1}{24\cdot 25}}, \quad\text{for each $k\geq 1$ and $h\geq 2$.}
\]
and that
\[
\left(1\pm\frac{1}{n}\right)^n \leq\mathrm e^{\pm 1},\quad n\geq 1.
\]
By this latter inequality (with the sign $+$) and the Stirling lower bound of the factorial,  we have 
\[
\frac{\mathrm{e}^{-k+1}}{(k-1)!} k^{k-1} \leq \frac{\mathrm e}{\sqrt {2\pi (k-1)}},\quad k\geq 2.
\]
Therefore,  for $k\geq 2$ and $h\geq 2$,  we have
\begin{equation}
\label{ineq:binom_coeff}
\binom{kh}{k-1}\leq\frac{\mathrm{e}^{\frac{1}{24\cdot 25}}}{h-1}\frac{1}{\sqrt{2\pi(k-1)}}\sqrt{1+\frac{1}{h-1}}
\left(\frac{h^h}{(h-1)^{h-1}}\right)^{k}=:C_{h,k}.
\end{equation}
We continue the proof distinguishing two cases: $h=1$ and $h\geq 2$.\\
\noindent$\it{Case\,\,h=1.}$\\ 
Since
\[
\binom{kh}{k-1}=k
\]
by \eqref{eq:03112023IV},  for any $m\in\mathbb N$,  we have
\[
\mathbb E Z^m=\frac{1-p}{p}\sum_{k\geq 1} k^m p^{k}=\frac{1-p}{p}\sum_{k\geq 1} k^m\mathrm{e}^{-\nu_1 k},
\quad\text{where $\nu_1:=-\log p>0$.}
\]

Using Lemma \ref{lemma:abel-plana} we have
\begin{equation}\label{eq:03112023V}
\frac{\E Z^m}{\sqrt{\lambda_\ell \mathrm{Leb}(B_\ell \cap C_\ell)}^{m-2}} \leq \frac{(1-p)m!}{p\sqrt{\lambda_\ell \mathrm{Leb}(B_\ell \cap C_\ell)}^{m-2}}\left(\frac{1}{\nu_1^{m+1}}+\frac{1}{\pi m m!} + \frac{2}{\pi^{m+1}} \right),\quad m,\ell\in\mathbb N.
\end{equation}
\noindent{\it Proof\,\,of\,\,Part\,\,(i)1}.\\
If $p\leq\mathrm{e}^{-1}$,  then $\nu_1 \geq 1$,  therefore,  setting  
$u_p:=1.05 \frac{1-p}{p}$,
by \eqref{eq:03112023V}, for any $m\geq 3$ and $\ell\in\mathbb N$,  we have
\begin{align}
\frac{\E Z^m}{\sqrt{\lambda_\ell \mathrm{Leb}(B_\ell \cap C_\ell)}^{m-2}} &\leq \frac{m!}{\sqrt{\lambda_\ell \mathrm{Leb}(B_\ell \cap C_\ell)}^{m-2}}1.05 \frac{1-p}{p}\nonumber\\
&=\frac{m!}{\sqrt{\lambda_\ell \mathrm{Leb}(B_\ell \cap C_\ell)}^{m-2}}u_p\nonumber\\
&=\frac{m!}{((u_p^{-1})\sqrt{\lambda_\ell \mathrm{Leb}(B_\ell \cap C_\ell)})^{m-2}}(u_p^{-1})^{m-3}\nonumber\\
&\leq\frac{m!}{((u_p^{-1})\sqrt{\lambda_\ell \mathrm{Leb}(B_\ell \cap C_\ell)})^{m-2}},\label{eq:15122023pom1}
\end{align}
where the latter inequality follows noticing that $u_p^{-1}\leq 1$,  indeed
with $\frac{1}{p}-1 \geq \log \left( \frac{1}{p}\right) \geq 1$,  and so $u_p\geq 1$.
Combining \eqref{eq:15122023pom1} with \eqref{eq:newassump} we have that condition \eqref{hyp:29092023} is satisfied with 
$\Delta_\ell := \frac{p}{1.05(1-p)}\sqrt{\lambda_\ell \text{Leb}(B_\ell\cap C_\ell)}$,  and the claim follows by Corollary \ref{cor:ST}.\\
\noindent{\it Proof\,\,of\,\,Part\,\,(i)2.}\\
If $p>\mathrm{e}^{-1}$, then $\nu_1<1$,  therefore,  setting  
$u_p:=\frac{1-p}{p\nu_1^3}$,  by \eqref{eq:03112023V}, for any $m\geq 3$ and $\ell\in\mathbb N$,  
\begin{align*}
    \frac{\E Z^m}{\sqrt{\lambda_\ell \mathrm{Leb}(B_\ell \cap C_\ell)}^{m-2}} &\leq \frac{(1-p)m!}{p(\nu_1\sqrt{\lambda_\ell \mathrm{Leb}(B_\ell \cap C_\ell)})^{m-2}}\nu_1^{m-2}\left(\frac{1}{\nu_1^{m+1}}+\frac{1}{\pi m m!} + \frac{2}{\pi^{m+1}} \right)\\
    &= \frac{1-p}{p\nu_1^3} \frac{m!}{(\nu_1\sqrt{\lambda_\ell \mathrm{Leb}(B_\ell \cap C_\ell)})^{m-2}}\left(1+\frac{\nu_1^{m+1}}{\pi m m!} + 2\left(\frac{\nu_1}{\pi}\right)^{m+1} \right)\\
    &\leq \frac{1-p}{p\nu_1^3} \frac{m!}{(\nu_1\sqrt{\lambda_\ell \mathrm{Leb}(B_\ell \cap C_\ell)})^{m-2}}1.05\\
&=u_p\frac{m!}{(\nu_1\sqrt{\lambda_\ell \mathrm{Leb}(B_\ell \cap C_\ell)})^{m-2}}1.05\\
&=[(1.05 u_p)^{-1}]^{m-3}\frac{m!}{((1.05 u_p)^{-1}\nu_1\sqrt{\lambda_\ell \mathrm{Leb}(B_\ell \cap C_\ell)})^{m-2}}.
\end{align*}
Since the function $(\mathrm{e}^{-1},1)\ni p\mapsto u_p=\frac{1-p}{-p(\log p)^3}$ is increasing and $\lim_{p\to\mathrm{e}^{-1}}u_p=\mathrm{e}-1>1$ we have that $u_p>1$ and so
\[
1.05 u_p > 1,\quad\text{$p\in (\mathrm{e}^{-1},1)$.}
\]
Therefore, for all $m\geq 3$ and $\ell\in\mathbb N$, we have 
\begin{align*}
     \frac{\E Z^m}{\sqrt{\lambda_\ell \mathrm{Leb}(B_\ell \cap C_\ell)}^{m-2}} 
     &\leq  \frac{m!}{((1.05 u_p)^{-1}\nu_1\sqrt{\lambda_\ell \mathrm{Leb}(B_\ell \cap C_\ell)})^{m-2}}.
\end{align*}
Combining this latter inequality with \eqref{eq:newassump} we have that condition \eqref{hyp:29092023} is satisfied with 
$\Delta_\ell := \frac{p(\log p)^4}{1.05(1-p)}\sqrt{\lambda_\ell \text{Leb}(B_\ell \cap C_\ell)}$,  and the claim follows by Corollary \ref{cor:ST}.\\
\noindent$\it{Case\,\,h\geq 2.}$\\ 
If $h\geq 2$, then by \eqref{ineq:binom_coeff} we have
\[
\binom{kh}{k-1}\leq\ind_{\{k=1\}}+C_{h,k}\ind_{\{k\geq 2\}}.
\]
Combining this relation with \eqref{eq:03112023IV} we have

\begin{align*}
    \E Z^m &\leq (1-p)^h+\sqrt{1+\frac{1}{h-1}} \frac{\mathrm{e}^{\frac{1}{24\cdot 25}}}{(h-1)\sqrt {2\pi}}\sum_{k\geq 2} k^{m-1} \frac{1}{\sqrt {k-1}}\left(\frac{h^h}{(h-1)^{h-1}} \right)^kp^{k-1} (1-p)^{k(h-1)+1}\\
    &=(1-p)^h+\sqrt{1+\frac{1}{h-1}} \frac{\mathrm{e}^{\frac{1}{24\cdot 25}}(1-p)}{p(h-1)\sqrt {2\pi}}\sum_{k\geq 2} k^{m-1} \frac{1}{\sqrt {k-1}}\left(\frac{h^h}{(h-1)^{h-1}} \right)^kp^{k} (1-p)^{k(h-1)}\\
    &=(1-p)^h+\sqrt{1+\frac{1}{h-1}} \frac{\mathrm{e}^{\frac{1}{24\cdot 25}}(1-p)}{p(h-1)\sqrt {2\pi}}\sum_{k\geq 2} k^{m-1} \frac{1}{\sqrt {k-1}}\left(ph \left(\frac{h(1-p)}{h-1}\right) ^{h-1} \right)^k\nonumber\\
&\leq(1-p)^h+\sqrt{1+\frac{1}{h-1}} \frac{\mathrm{e}^{\frac{1}{24\cdot 25}}(1-p)}{p(h-1)\sqrt {2\pi}}\sum_{k\geq 2} k^{m-1}\mathrm{e}^{-\nu_2 k},
\quad\text{$m\geq 3$}
\nonumber
\end{align*}
where
\[
\nu_2:=-\log\left(ph\left( \frac{h(1-p)}{h-1}\right) ^{h-1}\right).
\]
Now we are going to verify that $\nu_2>0$, i.e.,  

\begin{equation}\label{eq:09102023due}
ph\left( \frac{h(1-p)}{h-1}\right) ^{h-1}<1.
\end{equation}
Setting $x:=ph \in (0,1)$,  we have
\begin{align*}
    x\left( \frac{h-x}{h-1}\right)^{h-1}&=x\left( \frac{h-1+1-x}{h-1}\right)^{h-1}\\
   &=x\left( 1+\frac{1-x}{h-1}\right)^{h-1}\\
   &= x\mathrm e^{(h-1)\log \left(1+ \frac{1-x}{h-1}\right)}\\
   &\leq x\mathrm e^{1-x}.
\end{align*}
Relation \eqref{eq:09102023due} follows noticing that the mapping
$x\in (0,1)\mapsto x\mathrm e^{1-x}$ is an increasing bijection from $(0,1)$ to itself.
Therefore, by Lemma \ref{lemma:abel-plana},  for any $m\geq 3$, we have
\begin{align}
\E Z^m 
&\leq(1-p)^h+(m-1)!\sqrt{1+\frac{1}{h-1}} \frac{\mathrm{e}^{\frac{1}{24\cdot 25}}(1-p)}{p(h-1)\sqrt {2\pi}}
\left(\nu_2^{-m}+\frac{1}{\pi (m-1)(m-1)!}+\frac{2}{\pi^m}\right)
\nonumber\\
&\leq(m-1)!\left((1-p)^h+\sqrt{1+\frac{1}{h-1}} \frac{\mathrm{e}^{\frac{1}{24\cdot 25}}(1-p)}{p(h-1)\sqrt {2\pi}}\right)
\left(\nu_2^{-m}+\frac{1}{\pi (m-1)(m-1)!}+\frac{2}{\pi^m}\right)\nonumber\\
&\leq(m-1)!u_{p,h}
\left(\nu_2^{-m}+\frac{1}{\pi (m-1)(m-1)!}+\frac{2}{\pi^m}\right),\label{eq:06112023primo}
\end{align}
where
\begin{equation*}
u_{p,h}:=1+\sqrt{1+\frac{1}{h-1}} \frac{\mathrm{e}^{\frac{1}{24\cdot 25}}(1-p)}{p(h-1)\sqrt {2\pi}}.
\end{equation*}
\noindent{\it Proof\,\,of\,\,Part\,\,(ii)1}.
If $ph(h(1-p)/(h-1))^{h-1}\leq\mathrm{e}^{-1}$,  i.e.,  $\nu_2\geq 1$, then combining \eqref{eq:13102023primo} (with $\nu_2$ in place of $\nu$) with \eqref{eq:06112023primo},  we have

\begin{align}
\frac{\mathbb E Z^m}{\sqrt{\lambda_\ell \text{Leb}(B_\ell\cap C_\ell)}^{m-2} }&\leq\frac{((u_{p,h})^{-1})^{m-3}m!}{((u_{p,h})^{-1}\sqrt{\lambda_\ell \text{Leb}(B_\ell \cap C_\ell)})^{m-2}}\nonumber\\
\leq & \frac{m!}{((u_{p,h})^{-1}\sqrt{\lambda_\ell \text{Leb}(B_\ell \cap C_\ell)})^{m-2}},\quad m\geq 3\label{eq:15122023pom2}
\end{align}
where we used that $(u_{p,h})^{-1}<1$.
Combining \eqref{eq:15122023pom2} with \eqref{eq:newassump} we have that condition \eqref{hyp:29092023} is satisfied with $\Delta_\ell:=(u_{p,h})^{-1}\sqrt{\lambda_\ell \text{Leb}(B_\ell \cap C_\ell)}$,  and the claim follows by Corollary \ref{cor:ST}.\\
\noindent{\it Proof\,\,of\,\,Part\,\,(ii)2}.  If $ph(h(1-p)/(h-1))^{h-1}>\mathrm{e}^{-1}$,  i.e.,  $\nu_2<1$,  then combining \eqref{eq:13102023secondo} (with $\nu_2$ in place of $\nu$) with \eqref{eq:06112023primo},
we have
\begin{align*}
\frac{\mathbb E Z^m}{\sqrt{\lambda_\ell \text{Leb}(B_\ell \cap C_\ell)}^{m-2}}&
    \leq u_{p,h}\frac{(m-1)!}{\sqrt{\lambda_\ell \text{Leb}(B_\ell \cap C_\ell)}^{m-2}}\left(\nu_2^{-m}+0.16\right)\nonumber\\
&=\frac{(m-1)!}{(\sqrt{\lambda_\ell \text{Leb}(B_\ell \cap C_\ell)}\nu_2)^{m-2}}\frac{\nu_2^m}{(u_{p,h})^{-1}\nu_2^2}\left(\nu_2^{-m}+0.16\right),\quad m\geq 3.\nonumber
\end{align*}

Since $\nu_2 <1$,  we have 
\begin{align*}
 \frac{\nu_2^m}{(u_{p,h})^{-1}\nu_2^2}\left(\nu_2^{-m}+0.16\right)\nonumber&\leq\frac{1.16}{(u_{p,h})^{-1}\nu_2^2}:=\widetilde{u}_{p,h}.
\end{align*}
Therefore
\begin{align*}
\frac{\mathbb E Z^m}{\sqrt{\lambda_\ell \text{Leb}(B_\ell \cap C_\ell)}^{m-2}}
&\leq\frac{(m-1)!}{(\sqrt{\lambda_\ell \text{Leb}(B_\ell \cap C_\ell)}\nu_2)^{m-2}}\widetilde{u}_{p,h}\nonumber\\
&=\frac{(m-1)!}{(\sqrt{\lambda_\ell \text{Leb}(B_\ell \cap C_\ell)}\nu_2 (\widetilde{u}_{p,h})^{-1})^{m-2}}((\widetilde{u}_{p,h})^{-1})^{m-3},\quad m\geq 3.
\end{align*}
Using Lemma \ref{lemma:inequality}, we have that $(\widetilde{u}_{p,h})^{-1} <(u_{p,h})^{-1}\nu_2^2<(u_{p,h})^{-1}<1$. Therefore
\begin{align*}
\frac{\mathbb E Z^m}{\sqrt{\lambda_\ell \text{Leb}(B_\ell \cap C_\ell)}^{m-2}} 
    \leq & \frac{m!}{(\sqrt{\lambda_\ell \text{Leb}(B_\ell \cap C_\ell)}\nu_2(\widetilde{u}_{p,h})^{-1})^{m-2}},
\quad m\geq 3.
\end{align*}
Combining this latter inequality with \eqref{eq:newassump} we have that condition \eqref{hyp:29092023} is satisfied with 
$\Delta_\ell :=\sqrt{\lambda_\ell \text{Leb}(B_\ell \cap C_\ell)}\nu_2(\widetilde{u}_{p,h})^{-1}$,  and the claim follows by Corollary \ref{cor:ST}.
\\
\noindent$\square$

\section{Application to a class of interferences in a wireless communication model}\label{sec:WIRE}

\subsection{Gaussian approximation}\label{sec:1512partic3}

In this section we apply Theorem \ref{thm:GaussPoissoncluster} to the interference $I(\{\bold 0\})$ (see e.g.  Remark \ref{re:particularcases}) when the Poisson process of nodes' locations has
a piecewise constant intensity function of the form $\lambda(x):=\lambda\bold{1}_B(x)$,  for some $\lambda>0$ and $B\in\mathcal{B}(\mathbb R^2)$.
In such a case we have quite explicit upper bounds on the Wasserstein and the Kolmogorov distances. 
The following corollary (whose proof is straightforward,  and therefore omitted) allows for explicit bounds 
for some classes of signal power distributions and attenuation functions.

\begin{Corollary}\label{thm:GaussInterference}
Let $\lambda(x):=\lambda\bold{1}_B(x)$,  $x\in\mathbb R^2$,  for some $\lambda>0$ and $B\in \mathcal{B}(\mathbb R^2)$ such that 
\[
0<\mathbb E Z_1^2\int_B A^2(x)\mathrm{d}x<\infty.
\]
Then
\[
d_W(L\{\bold 0\}),G)\leq\frac{1}{\sqrt\lambda}\frac{\mathbb E Z_1^3}{(\mathbb E Z_1^2)^{3/2}}\frac{\int_B A(x)^3\mathrm{d}x}{\left(\int_B A(x)^2\mathrm{d}x\right)^{3/2}}
\]
and
\begin{align}
&d_K(L(\{\bold 0\}),G)\nonumber\\
&\qquad
\leq\frac{1}{\sqrt\lambda}
\left[1+\frac{1}{2}\max\Biggl\{4,\left[4\frac{1}{\lambda}\frac{\mathbb E Z_1^4}{(\mathbb E Z_1^2)^2}\frac{\int_{B}A(x)^4\mathrm{d}x}{\left(\int_{B}A(x)^2\mathrm{d}x\right)^2}+2\right]^{1/4}\Bigg\}\right]\frac{\mathbb E Z_1^3}{(\mathbb E Z_1^2)^{3/2}}\frac{\int_{B}A(x)^3\mathrm{d}x}{\left(\int_{B}A(x)^2\mathrm{d}x\right)^{3/2}}\nonumber\\
&\qquad\qquad\qquad
+\frac{1}{\sqrt\lambda}\frac{\sqrt{\mathbb E Z_1^4}}{\mathbb E Z_1^2}\frac{\sqrt{\int_{B}A(x)^4\mathrm{d}x}}{\int_{B}A(x)^2\mathrm{d}x}.
\nonumber
\end{align}
\end{Corollary}

\begin{Example}
If the path loss function is the Hertzian attenuation function, i.e.,  $A(x):=\max\{R,\|x\|\}^{-\alpha}$,  $x\in\R^2$,  for some $R>0$ and $\alpha>1$,  $B:=\R^2$
and $\E Z_1^2\in (0,\infty)$,  then Corollary \ref{thm:GaussInterference} applies with
\[
\int_B A(x)^m\mathrm{d}x=2\pi\left(\frac{1}{2}-\frac{1}{2-\alpha m}\right)R^{2-\alpha m},\quad\text{for $m=2,3,4$.}
\]
\end{Example}

\subsubsection{Moderate deviations,  Bernstein-type concentration inequalities and normal approximation bounds with Cram\'er correction term}\label{sec:1512partic4}

In this section we apply Theorem \ref{thm:ModeratePoissoncluster} to the sequence $\{I_\ell(\{\bold 0\})\}_{\ell\geq 1}$ (defined in Remark \ref{re:particularcasesbis})
when the Poisson processes of nodes' locations have a piecewise deterministic intensity function of the form $\lambda_\ell(x):=\lambda_\ell\bold{1}_{B_\ell}(x)$,  $x\in\mathbb R^2$, for some sequences $\{\lambda_\ell\}_{\ell\geq 1}\subset (0,\infty)$ and $\{B_\ell\}_{\ell\geq 1}\subset\mathcal{B}(\mathbb R^2)$.
In such a case the assumption \eqref{eq:08122023primo} simplifies a lot. 

The following corollary (whose proof is straightforward, and therefore omitted) holds.

\begin{Corollary}\label{cor:ModerateInterference}
Let $\{B_\ell\}_{\ell \in \mathbb N}\subset\mathcal{B}(\mathbb R^d)$ and $\{\mathbb Q_\ell\}_{\ell\geq 1}$ be such that 
\begin{equation}\label{eq:19122023I}
0<\mathbb E (Z_1^{(\ell)})^2\int_{B_\ell}A(x)^2\mathrm{d}x<\infty,\quad \ell\geq 1
\end{equation}
and assume that there exist a non-negative constant $\gamma\geq 0$ and a positive numerical sequence $\{\Delta_\ell\}_{\ell \in \mathbb N}$ such that 
\begin{equation}
\label{eq:19122023II}
\frac{1}{\lambda_\ell^{\frac{m}{2}-1}}\frac{\mathbb E(Z_1^{(\ell)})^m}{(\mathbb{E}(Z_1^{(\ell)})^2)^{m/2}}
\frac{\int_{B_\ell}A(x)^m\mathrm{d}x}{\left(\int_{B_\ell}A(x)^2\mathrm{d}x\right)^{m/2}}
\leq \frac{(m!)^{1+\gamma}}{\Delta_\ell ^{m-2}},\quad \text{for all } m\geq 3 \text{ and }\ell\in\mathbb N.
\end{equation}
Then the sequence $\{I_\ell(\{\bold 0\})\}_{\ell\geq 1}$ satisfies a $\bold{MDP}(\gamma,\{\Delta_\ell\}_{\ell\in\mathbb N})$,  a $\bold{BCI}(\gamma,\{\Delta_\ell\}_{\ell\in\mathbb N})$ and a\\ $\bold{NACC}(\gamma,\{\Delta_\ell\}_{\ell\in\mathbb N})$.
\end{Corollary}

\begin{Example}
Under the notation of Corollary \ref{cor:ModerateInterference},  if we set $B_\ell\equiv\mathbb R^2$ for any $\ell\geq 1$,  we suppose that $\mathbb Q_\ell$ is the exponential law with mean $\mu^{-1}$, for some $\mu>0$, and we assume that the attenuation of the signal is Hertzian, i.e.,  $A(x):=\max\{R,\|x\|\}^{-\alpha}$,  $x\in\mathbb R^2$, for some constants $R>0$ and $\alpha>1$,  then
\[
\frac{\mathbb E(Z_1^{(\ell)})^m}{(\mathbb{E}(Z_1^{(\ell)})^2)^{m/2}}=\frac{m!}{2^{m/2}}\quad\text{and}\quad\int_{\R^2}A(x)^m\mathrm{d}x=\frac{\pi\alpha m}{\alpha m-2}R^{2-\alpha m},\quad\text{for any $\ell\in\mathbb N$ and $m\geq 2$.}
\]
So the assumption \eqref{eq:19122023I} of Corollary \ref{cor:ModerateInterference} is satisfied and the left-hand side of relation \eqref{eq:19122023II} reads, for any $\ell\in\mathbb N$ and $m\geq 3$,
\[
\frac{m!}{\left(R\sqrt{2\lambda_\ell\frac{\pi\alpha}{\alpha-1}}\right)^{m-2}}\times\frac{m(2\alpha-2)}{4(\alpha m-2)}.
\]
Since $m\geq 3$ we have
\[
\frac{\alpha m-m}{\alpha m-2}<1
\]
So Corollary \ref{cor:ModerateInterference} yields that the sequence $\{I_\ell(\{\bold 0\})\}_{\ell\geq 1}$ satisfies a $\bold{MDP}(0,\{\Delta_\ell\}_{\ell\in\mathbb N})$,  a $\bold{BCI}(0,\{\Delta_\ell\}_{\ell\in\mathbb N})$ and a $\bold{NACC}(0,\{\Delta_\ell\}_{\ell\in\mathbb N})$ with
\[
\Delta_\ell:=R\sqrt{2\lambda_\ell\frac{\pi\alpha}{\alpha-1}},\quad\ell\geq 1.
\]
\end{Example}
\section{Conclusion}
Exploiting the theory developed in \cite{LPS}, we provided explicit bounds on the Wasserstein and the Kolmogorov distances
between random variables lying in the first chaos of the Poisson space and the standard Normal
distribution. Relying on the findings in \cite{SS} and on a
fine control of the cumulants of the first chaos on the Poisson space, we also provided moderate deviations, Bernstein-type
concentration inequalities and Normal approximation bounds with Cramér correction terms for
the same random variables.  We applied these results to Poisson shot-noise random variables,
and in particular to generalized compound Hawkes point processes. As far as Hawkes processes is concerned,
%If we focus on Hawkes processes on $(0,+\infty)$,  
the results proven in this paper generalise many of the asymptotic theorems found in the literature \cite{HHKR, KPR, GZ,ZHU} to the spatial case, eventually with a varying baseline intensity and with less constraining assumptions on the excitation kernels.

We point out that some Hawkes processes have a Galton-Watson representation but cannot be easily expressed as a Poisson integral of the type \eqref{eq:Pois_integral}. The main example is that of a multivariate Hawkes process exhibiting self and cross excitation between many interacting nodes. Indeed, such a process does have a branching structure \cite{ELL} but \textit{a priori} does not fall within the context of this paper.  To the best of our knowledge,  we only have bounds on the Wasserstein metric between multivariate Hawkes processes with exponential kernels and their multivariate Gaussian limit which is of order $O\left ( 1/\sqrt t \right)$ \cite{Khabou}.

Another interesting development of the results proven in this paper would be their extension to the whole path of the process, rather than the process evaluated at one instant. More specifically, we would like to find upper bounds on the distance between the centered and normalised \textit{path} of the Poisson shot noise process, and its limiting Gaussian \textit{process} in the space of càdlàg functions equipped with the Skorokhod metric, for example using the results provided in \cite{BRZ}. 
These approximation results are obviously more delicate to obtain, and to the best of our knowledge they have only been studied in few works such as \cite{BCDM}. 

\section{Proofs of Lemmas \ref{lemma:Z},  \ref{lemma:abel-plana} and \ref{lemma:inequality}}\label{sec:lemmas}

\subsection{Proof of Lemma \ref{lemma:Z}}

The claim is clearly true if $\max\{\mathrm{Leb}(B\cap C),\mathbb E Z^m,\E M^m\}=+\infty$.
Therefore we assume  $\max\{\mathrm{Leb}(B\cap C),\E Z^m,\E M^m\}<+\infty$. We start with the obvious inequality
\begin{equation}\label{eq:14122023primo}
\E |v(Z_1)(C-x)|^m\leq\E \left(\sum_{k=0}^{Z_1(C-x,\R)-1}|M_{k,1}|\right)^m.
\end{equation}

Using Hölder's inequality we have 
\begin{align*}
    \sum_{k=0}^{Z_1(C-x,\R)-1}|M_{k,1}|&\leq \left(\sum_{k=0}^{Z_1(C-x,\R)-1} 1 \right)^{\frac{m-1}{m}} \left( \sum_{k=0}^{Z_1(C-x,\R)-1}|M_{k,1}|^m\right)^{\frac{1}{m}}
\end{align*}
which yields after raising the inequality to the $m-$th power 
$$\left(\sum_{k=0}^{Z_1(C-x,\R)-1}|M_{k,1}| \right)^m \leq \left(Z_1(C-x,\R) \right)^{m-1}  \sum_{k=0}^{Z_1(C-x,\R)-1}|M_{k,1}|^m.$$
Using the independence between $Z_1(C-x,\R)$ and $(|M_{k,1}|)_{k \in \mathbb N}$ and Wald's identity we have 
\begin{align*}
    \E \left(\sum_{k=0}^{Z_1(C-x,\R)-1}|M_{k,1}|\right)^m &\leq \E \left[ \left(Z_1(C-x,\R) \right)^{m-1}  \sum_{k=0}^{Z_1(C-x,\R)-1}|M_{k,1}|^m\right]\\
    &= \E \left[\left(Z_1(C-x,\R) \right)^{m-1} \E \left[   \sum_{k=0}^{Z_1(C-x,\R)-1}|M_{k,1}|^m \bigg | Z_1(C-x,\R)\right ]\right]\\
    &=\E \left[\left(Z_1(C-x,\R) \right)^{m} |M|^m\right],
\end{align*}
and finally inequality \eqref{eq:14122023primo} yields
$$\E |v(Z_1)(C-x)|^m\leq \E |M|^m \E Z_1(C-x,\R)^m.$$

Recalling that we denote by $\{Y_{1,k}\}_{k\geq 0}$,  $Y_{1,0}:=\bold 0$,  the first components of the points of $Z_1(\cdot,\cdot)$,  we have 

    \begin{align}
        Z_1(C-x,\R)^m&=\card \left \{k \in \mathbb N\cup\{0\}: Y_{1,k} \in C-x \right \}^m\nonumber \\
        &=\left(\sum_{k=0}^{Z-1} \boldsymbol 1_{Y_{1,k} \in C-x }\right)^m,\quad m\geq 1.\label{eq:23052023secondo}
    \end{align}
    The $m$th power of the sum of indicators can be expanded by using the multinomial theorem, which yields
    \begin{align}
        \left(\sum_{k=0}^{Z-1} \boldsymbol 1_{Y_{1,k} \in C-x }\right)^m &= \sum_{\substack{k_0,\ldots,k_{Z-1} \geq 0 \\
        k_0+\cdots+ k_{Z-1}=m}} \binom{m}{k_0,\ldots, k_{Z-1}} \prod _{i=0}^{Z-1} \boldsymbol 1_{Y_{1,i} \in C-x }^{k_i}\nonumber \\
        &= \sum_{\substack{k_0,\ldots,k_{Z-1} \geq 0 \\
        k_0+\cdots+ k_{Z-1}=m}} \binom{m}{k_0,\ldots, k_{Z-1}} \prod _{i=0}^{Z-1} \boldsymbol 1_{Y_{1,i} \in C-x }\nonumber \\
        &\leq \sum_{\substack{k_0,\ldots,k_{Z-1} \geq 0 \\
        k_0+\cdots+ k_{Z-1}=m}} \binom{m}{k_0,\ldots, k_{Z-1}}  \boldsymbol 1_{Y_{1,0} \in C-x },\quad m\geq 1.\label{ineq:multinomial}
    \end{align}
Here
\[
\binom{m}{k_0,\ldots, k_{Z-1}}=\frac{m!}{k_0!\ldots k_{Z-1}!}
\]
denotes the multinomial coefficient. By \eqref{ineq:multinomial} we have

    \begin{align}
        \int_{B}  \left(\sum_{k=0}^{Z-1} \boldsymbol 1_{Y_{1,k} \in C-x }\right)^m  \mathrm d x &\leq \sum_{\substack{k_0,\ldots,k_{Z-1} \geq 0 \\
        k_0+\cdots+ k_{Z-1}=m}} \binom{m}{k_0,\ldots, k_{Z-1}}  \int_{B}\boldsymbol 1_{Y_{1,0} \in C-x } \mathrm d x \nonumber\\
&=\sum_{\substack{k_0,\ldots,k_{Z-1} \geq 0 \\
        k_0+\cdots+ k_{Z-1}=m}} \binom{m}{k_0,\ldots, k_{Z-1}}  \int_{B}\boldsymbol 1_{\bold 0 \in C-x } \mathrm d x \nonumber\\
        &= \sum_{\substack{k_0,\ldots,k_{Z-1} \geq 0 \\
        k_0+\cdots+ k_{Z-1}=m}} \binom{m}{k_0,\ldots, k_{Z-1}} \mathrm {Leb}(B\cap C) \nonumber \\
         &= \mathrm {Leb}(B\cap C) \sum_{\substack{k_0,\ldots,k_{Z-1} \geq 0 \\
        k_0+\cdots+ k_{Z-1}=m}} \binom{m}{k_0,\ldots, k_{Z-1}} \prod _{i=0}^{Z-1} 1^{k_i} \nonumber\\
&\leq\mathrm {Leb}(B\cap C) Z^m, \quad m\geq 1\label{ineq:third}
    \end{align}
where the latter inequality follows by using again the multinomial theorem.
The claim easily follows by \eqref{eq:14122023primo},  \eqref{eq:23052023secondo} and \eqref{ineq:third}.

\subsection{Proof of Lemma \ref{lemma:abel-plana}}

Set $D:=\{z\in\mathbb C:\,\,\mathrm{Re}z>0\}$ and define $f(z):=z^{m-1}\mathrm{e}^{-\nu z} $,  $z\in D$,  $m\geq 2$,  $\nu>0$.

Clearly,  $f$ is analytic on $D$; we shall check later on that 
\begin{equation}\label{eq:06112023I}
\text{For any compact $K\subset (0,\infty)$,  $\lim_{y\to +\infty}\sup_{x\in K}|f(x\pm \bold i y)|e^{-2\pi y}=0$}
\end{equation}
and
\begin{equation}\label{eq:06112023II}
\text{For any $x>0$,\,\,}\int_0^{+\infty} |f(x+\bold i y)-f(x-\bold i y)|\mathrm{e}^{-2\pi y } \mathrm d y<\infty\quad\text{and}\quad\lim_{x\to+\infty}\int_0^{+\infty} |f(x+\bold i y)-f(x-\bold i y)|\mathrm{e}^{-2\pi y } \mathrm d y=0.
\end{equation}
Therefore by the Abel-Plana formula (see e.g.  \cite{BFSS}) we have (note that $f(0)=0$)

\begin{align}
\sum_{k\geq 1}f(k)&=\int_0^\infty f(t)\d t+\bold i \int_0^\infty\frac{f(\bold i t)-f(-\bold i t)}{\mathrm e^{2\pi t}-1}\d t\nonumber\\
&= \int_0^{+\infty}  \mathrm e^{-\nu t} t^{m-1} \dM t + \bold i \int_0^{+\infty} \frac{\mathrm e^{-\bold i\nu t}(\bold i t)^{m-1}-\mathrm e^{\bold i \nu t}(- \bold i t)^{m-1}}{\mathrm e^{2\pi t}-1} \dM t\nonumber\\
    &= \nu^{-m} (m-1)!+ R_m,\label{eq:06112023quinto}
\end{align}
where we used that
\begin{equation}\label{eq:06112023X}
\frac{\nu^m}{\Gamma(m)}\int_0^\infty t^{m-1}\mathrm{e}^{-\nu t}\mathrm{d}t=1
\end{equation}
and that the Euler gamma function $\Gamma(\cdot)$ computed on the integer $m$ is equal to $(m-1)!$.

We proceed by bounding $|R_m|$ from above.  We distinguish two cases: $m=2p$ and $m=2p+1$,  $p\in \mathbb N$.  If $m=2p$ we have
\begin{align*}
    R_{2p}&= \bold i \int_0^{+\infty} \frac{\mathrm e^{-\bold i\nu t}(\bold it)^{2p-1}-\mathrm e^{\bold i \nu t}(-\bold it)^{2p-1}}{\mathrm e^{2\pi t}-1} \dM t\\
    &=\bold i \int_0^{+\infty} \frac{\mathrm e^{-\bold i\nu t}(\bold i t)^{2p-1}+\mathrm e^{\bold i\nu t}(\bold i t)^{2p-1}}{\mathrm e^{2\pi t}-1} \dM t\\
    &=\bold i^{2p} \int_0^{+\infty} \frac{ t^{2p-1} (\mathrm e^{-\bold i \nu t}+\mathrm e^{\bold i\nu t})}{\mathrm e^{2\pi t}-1} \dM t\\
    &=2(-1)^p \int_0^{+\infty} \frac{ t^{2p-1} \cos(\nu t)}{\mathrm e^{2\pi t}-1} \dM t.\\
\end{align*}
Thus
\begin{align}
    |R_m|& \leq 2\int_0^{+\infty} \frac{ t^{m-1}}{\mathrm e^{2\pi t}-1} \dM t\nonumber\\
    &= 2\int_0^{1} \frac{ t^{m-1}}{\mathrm e^{2\pi t}-1} \dM t +2 \int_1^{+\infty} \frac{ t^{m-1}}{\mathrm e^{2\pi t}-1} \dM t\nonumber\\
    &\leq \frac{1}{  \pi}\int_0^{1} \frac{ t^{m-1}}{t} \dM t + 2\int_1^{+\infty} \mathrm e^{-\pi t} t^{m-1} \dM t\nonumber\\ 
    &\leq \frac{1}{\pi(m-1)} + 2\pi^{-m} (m-1)!,\quad m=2p\label{eq:06112023sesto}
\end{align}
where the latter inequality follows by \eqref{eq:06112023X} with $\pi$ in place of $\nu$.  Similarly,  we have 
\begin{align}
    |R_{m}| &= \left | \bold i \int_0^{+\infty} \frac{\mathrm e^{-\bold i\nu t}(\bold i t)^{m-1}-\mathrm e^{\bold i\nu t}(-\bold i t)^{m-1}}{\mathrm e^{2\pi t}-1} \dM t \right |\nonumber\\
    &= 2 \left |\int_0^{+\infty} \frac{ t^{m-1} \sin(\nu t)}{\mathrm e^{2\pi t}-1} \dM t\right |
\nonumber\\
&\leq 2 \int_0^{+\infty} \frac{ t^{m-1}}{\mathrm e^{2\pi t}-1} \dM t
\nonumber\\
    &\leq \frac{1}{\pi(m-1)} + 2\pi^{-m} (m-1)!,\quad m=2p+1.\label{eq:06112023settimo}
\end{align}
The claim follows by the relations \eqref{eq:06112023quinto},  \eqref{eq:06112023sesto} and \eqref{eq:06112023settimo}.

It remains to prove \eqref{eq:06112023I} and \eqref{eq:06112023II}.  We start proving \eqref{eq:06112023I}.  Let $K\subset (0,\infty)$ be an arbitrary compact set.  We have  
    \begin{align*}
        \sup_{x\in K}|f(x\pm \bold i y)|\mathrm{e}^{-2\pi y}&=\mathrm{e}^{-2\pi y}\sup_{x\in K}(x^2+y^2)^{\frac{m-1}{2}}\mathrm{e}^{-\nu x} \\
        &\leq ((\sup K)^2+y^2)^{\frac{m-1}{2}}\mathrm{e}^{-2\pi y}\to 0,\quad\text{as $y\to+\infty$.}
    \end{align*}
Finally we prove \eqref{eq:06112023II}.  For any $x,y>0$,  we have
    \begin{align*}
        |f(x+\bold i y)-f(x-\bold i y)| \mathrm{e}^{-2\pi y } &=|(x+\bold i y)^{m-1}\mathrm e^{-\nu\bold i y}-(x-\bold i y)^{m-1}\mathrm e^{\nu\bold i y}|e^{-2\pi y }\mathrm e^{-\nu x}\\
        &=\left |\sum_{k=0}^{m-1} \binom{m-1}{k}x^k \left((\bold i y)^{m-1-k}\mathrm e^{-\nu\bold i y} -(-\bold i y)^{m-1-k}\mathrm e^{\nu\bold i y}  \right) \right |\mathrm e^{-2\pi y } \mathrm e^{-\nu x}\\
        &\leq 2 \sum_{k=0}^{m-1} \binom{m-1}{k}x^k \mathrm e^{-\nu x} y^{m-1-k} \mathrm e^{-2\pi y }.
    \end{align*}
    Therefore,  for any $x>0$, we have
\begin{align}
\int_0^{+\infty} |f(x+\bold i y)-f(x-\bold i y)| \mathrm e^{-2\pi y } \mathrm d y &\leq 2 \sum_{k=0}^{m-1} \binom{m-1}{k}x^k \mathrm e^{-\nu x} \int_0^{+\infty}y^{m-k-1} \mathrm e^{-2\pi y }\mathrm d y\nonumber\\
&=  2 \sum_{k=0}^{m-1} \binom{m-1}{k}x^k \mathrm e^{-\nu x}\frac{(m-k-1)!}{(2\pi)^{m-k}},\label{eq:06112023XII}
\end{align}
where the latter equality follows by the relation \eqref{eq:06112023X} with $m-k$ in place of $m$ and $2\pi$ in place of $\nu$. Clearly,  the right-hand side of the relation
\eqref{eq:06112023XII} is finite and tends to zero as $x\to+\infty$.  The proof is completed.

\subsection{Proof of Lemma \ref{lemma:inequality}}

A simple computation shows
\[    
f'(x)= (x-1-\log x)\left(3x-3-\log x\right),\quad x\in (0,1).
\]
Since $x-1-\log x>0$ for every $x\in (0,1)$,  the sign of $f'$ coincides with the sign of $g(x):=3x-3-\log x$,  $x\in (0,1)$.  Studying the derivative of $g$ we get that $g$ increases on $(1/3,1)$ and decreases on $(0,1/3)$ with a minimum at $x=1/3$.  Since
$\lim_{x\to 0^+}g(x)=+\infty$,  $g(1/3)<0$ and $\lim_{x\to 1^-}g(x)=0$,  we then have that there exists a unique $x^*\in (0,1/3)$ with $g(x^*)=0$ and $g$ is positive on $(0,x^*)$ and it is negative on $(x^*,1)$. Therefore $x^*$ is point of maximum of $f$, and consequently,  for any $x\in (0,1)$, we have

    \begin{align*}
        f(x) \leq f(x^*)
        = x^*(x^*-1-\log x^*)^2=4 x^*(1-x^*)^2\leq4 x^*(1-x^*)
        &\leq 8/9<1,
    \end{align*}
where we used that $g(x^*)=0$,  that the mapping $(0,1)\in x\mapsto x(1-x)$ increases on $(0,1/2)$ and that $x^*<1/3$.

\section*{Acknowledgement}
We wish to thank Prof. Matthias Schulte for his useful suggestions.\\
Mahmoud Khabou was supported by the Project EDDA (ANR-20-IADJ-0003) of the French National Research Agency (ANR); Giovanni Luca Torrisi was supported by group GNAMPA of INdAM. For the purpose of open access, the authors have applied a Creative Commons Attribution (CC BY) licence to any Author Accepted Manuscript version arising.

\end{document}